\documentclass[12pt,preprint]{elsarticle}
\usepackage[utf8]{luainputenc}
\usepackage{color}
\usepackage{array}
\usepackage{float}
\usepackage{mathtools}
\usepackage{multirow}
\usepackage{amsthm}
\usepackage{amsbsy}
\usepackage{amstext}
\usepackage{amssymb}
\usepackage{graphicx}
\usepackage{esint}

\makeatletter

\journal{ }

\providecommand{\tabularnewline}{\\}

\theoremstyle{plain}

\ifx\proof\undefined
\newenvironment{proof}[1][\protect\proofname]{\par
\normalfont\topsep6\p@\@plus6\p@\relax
\trivlist
\itemindent\parindent
\item[\hskip\labelsep\scshape #1]\ignorespaces
}{%
\endtrivlist\@endpefalse
}
\providecommand{\proofname}{Proof}
\fi

\@ifundefined{showcaptionsetup}{}{%
 \PassOptionsToPackage{caption=false}{subfig}}
\usepackage{subfig}
\makeatother

\providecommand{\theoremname}{Theorem}

\begin{document}
\begin{frontmatter}
\title{A remark concerning  divergence accuracy order for $\mathbf{H}(\mbox{div})$-conforming   finite element flux approximations}
       
\author[fec]{Philippe R. B. Devloo }

\ead{phil@fec.unicamp.br}

\author[Salinas]{Agnaldo M. Farias \corref{cor1}}

\ead{agnaldo.farias@ifnmg.edu.br}

\author[ime]{Sônia M. Gomes}

\ead{soniag@ime.unicamp.br}

\address[fec]{FEC-Universidade Estadual de Campinas, Campinas-SP, Brazil}

\address[Salinas]{Departamento de Matem\'{a}tica, IFNMG, Salinas,  MG, Brazil}

\address[ime]{IMECC-Universidade Estadual de Campinas, Campinas, SP, Brazil}

\cortext[cor1]{Corresponding author - Phone +55 19 982098044}

\begin{abstract}{\footnotesize The construction of finite element approximations in $\mathbf{H}(\mbox{div}, \Omega)$  usually requires the Piola transformation to  map vector polynomials   from a master element to vector fields in the elements  of a partition of the region $\Omega$.  It is known that  degradation may occur in convergence  order  if non affine geometric mappings are used.  On this point,  we revisit a general procedure for the improvement of two-dimensional  flux approximations discussed in a recent paper of this journal (Comput.  Math. Appl. 74 (2017) 3283–3295).  The starting point is an approximation scheme, which is known to provide  $L^2$-errors with accuracy of order $k+1$ for sufficiently smooth flux functions, and of order $r+1$ for  flux divergence.  An example is $RT_k$ spaces on quadrilateral meshes,  where $r=k$ or $k-1$ if linear or bilinear geometric isomorphisms are applied. Furthermore, the original space is required to be expressed by a factorization in terms  of  edge and internal shape flux functions. The goal is to  define a hierarchy of enriched flux approximations to reach arbitrary higher orders of  divergence accuracy  $r+n+1$ as desired, for any $n \geq 1$.  The enriched versions are defined by adding   higher degree internal shape functions of the original family of spaces at level $k+n$,  while keeping the original border fluxes at level $k$.  The case $n=1$  has been discussed in the mentioned publication for two particular examples.  General  stronger enrichment  $n>1$ shall  be analyzed and applied to Darcy's flow simulations, the global condensed systems to be solved  having same dimension and structure of the original scheme.}

\end{abstract}

\end{frontmatter}
\section{Introduction}
The paper is concerned with the construction of finite element approximation spaces in $
\mathbf{H}(\mbox{div}, \Omega)$ and their applications to flux approximations in mixed 
finite element methods for Darcy's flow problems \cite{BrezziFortin1991}.  Usually, the 
Piola transformation is used to  map vector polynomials defined in  master element to vector 
fields in the elements  of a partition of the region $\Omega$.  

It is known that  
degradation may occur in convergence  order  if non affine geometric mappings are used, as 
emphasized in \cite{ABF}, 
a publication devoted to the determination of optimal conditions 
required by the vector polynomials spaces on the master element in order to obtain optimal 
orders of accuracy on the flux and flux divergence on the mapped elements. Based on these 
findings, $ABF_k$ spaces were derived as enriched version of the classic Raviart-Thomas 
spaces $RT_k$  for general quadrilateral meshes \cite{RT}.   For rectangular meshes, both space configurations give $L^2$-errors  of order $k+1$ for flux and flux divergence. Using bi-linear mappings,  divergence 
 are reduced to order $k$ for $RT_k$ spaces,   but the convergence  order 
$k+1$ is restored when  $ABF_k$ spaces are applied. 

In this contribution,  we revisit a  procedure for the construction of improved  
two-dimensional  
flux approximations discussed in a recent paper of this journal \cite{Hdiv2D2017}, where  
enriched spaces  $RT_k^+$, for quadrilaterals, and $BDM_k^+$ (i.e., $BDFM_{k+1}$) and 
$BDM_k^{++}$, for triangles,    have been discussed. The goal is to extend the methodology  
to any given original approximation scheme providing  $L^2$-errors with accuracy of order $k
+1$ for sufficiently smooth flux functions, and of order $r+1$ for  flux divergence.  A 
hierarchy of enriched flux approximations is then defined in order to reach arbitrary higher 
orders of  divergence accuracy  $r+n+1$ as desired, for any $n \geq 1$.

The original space is required to be expressed by a factorization in terms  of  edge and 
internal shape functions. Then the enriched versions are constructed by adding   higher 
degree internal shape functions of the  original space at level $k+n$,  while keeping the 
original border fluxes at level $k$.  Thus, for sufficiently smooth vector functions, approximations based on  
 these enriched spaces maintain the original order $k+1$, and  their divergence can reach accuracy of higher order $r+n+1$, as $n$ increases. 

Consequently, when applied to  Darcy's flow simulations, stable convergent approximations are obtained with  order of accuracy  $k+1$ for the flux variable, and  accuracy of order $q+1$ for the potential variable, where  $q=\min \{k+1,r+n,t+n\}$, where $t$ is determined by the total degree of  polynomials included in the potential scalar space.  Typically, $q=k+1$ and the best order of accuracy for the potential is $k+2$. However,  divergence errors can get arbitrary higher orders $r+n+1$ as $n$ increases. It should be observed that the computational effort for matrix assembly is higher for the enriched frameworks, however their global condensed systems to be solved  have same dimension (and structure) of the original scheme, which is proportional to the space
dimension of the border fluxes for each element geometry.  

From the  implementation point of view, the main requirement is the 
possibility to establish the order of the edge flux functions independently. 
Once a higher order scheme has been implemented, the enriched spaces   are simply obtained 
by pruning  the shape functions at level $k+n$ whose  edge functions are of degree $>k$.  
As in the previous paper \cite{Hdiv2D2017}, taking  $BDM_{k}$ \cite{BDM} for triangles  and 
$RT_k$ \cite{RT} for quadrilaterals as original spaces, the implementation of their higher 
order enriched versions is implemented in the object-oriented 
scientific computational environment NeoPZ \cite{NeoPZ}.
In this implementation  vector shape functions are classified by edge or  internal types
 (see  \cite{Siqueira,Surface}). In the NeoPZ coding structure, it is possible to 
 choose the order of approximation of the edge flux functions independently from the order 
 of internal functions, a required  capability for the implementation of the enriched  flux 
 spaces  documented in this paper

The text is organized as follows. After setting the notation and general aspects concerning 
divergence conforming approximation spaces in Section \ref{section2},  the general space 
enrichment procedure is described in Section \ref{section3},  and a projection error 
analysis is performed in order to demonstrate  the enhancement of divergence accuracy given 
by the new spaces. Section \ref{section4} is dedicated to the application  of the proposed 
enriched approximations to the discretization of a mixed finite element formulation, and the 
results of numerical simulations confirm the theoretical error estimations derived in 
previous sections.  Section \ref{section5} gives the final conclusions of the article. 

\section{Notation and general aspects}\label{section2}

Let $\Omega\subset \mathbb{R}^2$,  be a computational region covered by a regular partition $\mathcal{T}_h=\{K\}$. For each geometric element $K \in \mathcal{T}_h$, there is an associated   master element
 $\hat{K}$ and an invertible geometric  diffeomorfism $F_{K}:\hat{K}\rightarrow K$ 
 transforming $\hat{K}$ onto $K$. For the present study, the master element is the triangle $\hat{K}= {\cal T} = \{ (\hat{x},\hat{y}); \hat{x}\geq 0, \hat{y} \geq 0, \hat{x}+\hat{y}\leq 1\}$, or the square  $\hat{K}= {\cal R}= [-1,1]\times[-1,1]$, and
$F_{K}$ is linear $F_{K}(\hat{x},\hat{y})=\mathbf{A}_{0}+\mathbf{A}_{1}\hat{x}+ \mathbf{A}_{2}\hat{y}$, or bi-linear  $F_{K}(\hat{x},\hat{y})= \mathbf{A}_{0}+ \mathbf{A}_{1}\hat{x}+ \mathbf{A}_{2}\hat{y}+ \mathbf{A}_3\hat{x}\hat{y}$, where  $\mathbf{A}_i$ are constant vectors in $\mathbb{R}^2$. 

For the examples considered in the present study, the   spaces defined on the master elements are defined in terms of polynomial of the following types: $\mathbb{P}_{k}$,  scalar polynomial space of total
degree $k$ (used for triangular elements); $\mathbb{Q}_{l,m}$ scalar polynomial space of maximum
degree $l$ in $x$, and $m$ in $y$  (used for quadrilateral elements).

\subsection{Approximation spaces}
 Based on  the partitions $\mathcal{T}_h$, finite dimensional approximation spaces  $\mathbf {V}_h \subset\mathbf{H}(\mbox{div},\Omega)$ and $U_h\subset L^2(\Omega)$ are piecewise defined over the elements $K$. Different contexts shall be considered for the  two types of element geometry, all  sharing the following basic characteristics.
\begin{enumerate}
\item A vector polynomial space $\mathbf{\hat{V}}$ and a scalar polynomial space $\hat{U}$ are considered on the master element $\hat{K}$.
\item In all the cases,  the  divergence operator maps $\mathbf{\hat{V}}$ onto $\hat{U}$:
\begin{equation}
\nabla\cdot \hat{\mathbf{V}} = \hat{U}. \label{DRhamComplex}
\end{equation}
\item The spaces $\hat{\mathbf{V}}$ are spanned by a hierarchy of vector shape functions, which are organized into two classes: the shape functions of interior type, with vanishing normal components over all element  edges, and the shape functions associated to the element edges, otherwise. Thus, a direct factorization $\hat{\mathbf{V}}=\hat{\mathbf{V}}^\partial  \oplus \mathring{\hat{\mathbf{V}}}$, in terms  of  face and internal flux functions, naturally holds.
\item The functions  $\mathbf{q} \in \mathbf {V}_h$ and $\varphi \in U_h$ are piecewise defined:  $\mathbf{q}|_K \in \mathbf{V}(K)$ and $\varphi|_K\in U(K)$, by   locally backtracking the polynomial spaces  $\hat{\mathbf{V}}$ and $\hat{U}$.  More precisely, the usual mappings $\mathbb{F}_K$, induced by the geometric diffeomorfisms $F_K$, give
\[U(K)=\mathbb{F}_K \hat{U}=\left\{ \varphi|\; \varphi \circ F_{K}\in\hat{U}\right\}. \] 
The Piola transformation $\mathbb{F}^{\text{div}}_K$ is used to form vector functions in $K$: 
\[
\mathbf{V}(K)=\mathbb{F}_{K}^{div}\hat{\mathbf{V}}=\left\{ \mathbf{q}|\;{J}_{K}\;DF_{K}^{-1}\mathbf{q}\circ F_{K}\in\hat{\mathbf{V}}\right\} ,\label{DivPiola}
\]
where $DF_{K}$ is the Jacobian matrix of $F_{K}$, and ${J}_{K}=det(DF_{K})$. It is asssumed that $ {J}_{K}>0$. It can be verified that for $\hat{\mathbf{q}} \in \hat{\mathbf{V}}$, and  $\hat{\mathbf{x}} \in \hat{K}$, $\nabla \cdot \hat{\mathbf{q}}={J}_{K}(\hat{\mathbf{x}}) \nabla \cdot \mathbf{q}$.

\end{enumerate}

\subsubsection*{Remarks}
\begin{itemize}
\item The construction of hierarchic shape functions  $\hat{\boldsymbol{\Phi}}$  is described in \cite{Siqueira,Surface} for  some classic vector spaces $\hat{\mathbf{V}}$ based on   triangle ($BDM_k$ and $BDFM_k$) and quadrilateral ($RT_k$) master elements. It is based on  appropriate choice of constant vector fields $\hat{\mathbf{v}}$, based on the geometry of each master element, which are multiplied by an available set of $H^1$ hierarchical scalar basis functions $\hat{\varphi}$ to get $\hat{\boldsymbol{\Phi}}= \hat{\mathbf{v}}\hat{\varphi}$. There are shape functions of interior type, with vanishing normal components over all element faces. Otherwise, the shape functions are classified of edge type. The normal component of an edge function coincides on  the corresponding edge with the associated scalar shape function, and vanishes over the other edges. 
Particularly, the flux spaces   $\hat{\mathbf{V}}$ are spanned by  bases $ \hat{\mathbf{B}}$ of the form
\[
\begin{footnotesize}
\begin{array}[t]{ccccc}
\hat{\mathbf{B}} & = & \underbrace{\left\{ \boldsymbol{\Phi}^{\hat{l}_{m},\hat{a}_{s}},\;\boldsymbol{\Phi}^{\hat{l}_{m},n}\right\} } & \cup & \underbrace{\left\{ \boldsymbol{\Phi}^{\hat{K},\hat{l}_{m},n},\boldsymbol{\Phi}_{(1)}^{\hat{K},n_{1},n_{2}},\;\boldsymbol{\Phi}_{(2)}^{\hat{K},n_{1},n_{2}}\right\} }, \\
 &  & \mbox{edge functions} &  & \mbox{internal functions}
\end{array}
\end{footnotesize}
\]
where $\hat{l}_{m}$ are edges of $\hat{K}$, $\hat{a}_{s}$ are vertices of $\hat{l}_{m}$, $n$, $n_{1}$ and $n_{2}$ determine the degree of the shape functions, and the subscripts $(i), i=1,2$, indicate two linearly independent vector fields $\hat{\bf{v}}_{(i)}$ used in the definition of internal functions. This construction has been extended in \cite{Hdiv3D} for  tetrahedral, affine hexahedral and prismatic elements. See also \cite{Surface} for a detailed explanation and application to curvilinear two dimensional meshes on manifolds.

\item Therefore,  the direct factorization $\hat{\mathbf{V}} = \hat{\mathbf{V}}^{\partial} \oplus\mathring{\hat{\mathbf{V}}}$ 	
is naturally identified for these vector spaces, where
$\mathring{\hat{\mathbf{V}}}$ is the space   spanned by  the internal shape functions, and  $\hat{\mathbf{V}}^{\partial}$ being its complement, spanned by the edge shape functions.

\item Usually, the divergence space $\nabla \cdot \hat{\mathbf{V}}=\hat{U}=\hat{U}_0\oplus \hat{U}^\perp$, where $\hat{U}_0$ are the  constant functions, and $\hat{U}^\perp$ are the functions with zero mean, the image of the internal functions by the divergence operator, $\hat{U}^\perp = \nabla \cdot \mathring{\hat{\mathbf{V}}}$. 

\end{itemize}

\subsection{General approximation theory}
Classic error analyses for  finite element methods are usually based on estimations in terms of certain projection errors of the exact solution {\color{blue}onto} the approximation spaces.

\subsubsection*{The projections}

For the applications of the present study, bounded projections $\hat{\boldsymbol{\pi}}:\mathbf{H}^{\alpha}(\hat{K})\rightarrow\hat{\mathbf{V}}$ are considered for sufficiently smooth vector functions  $\hat{\mathbf{q}} \in \mathbf{H}^{\alpha}(\hat{K})$, $\alpha \geq 1$,  defined on the master element,  and they corresponding versions  $\boldsymbol{\pi}_{K}:\mathbf{H}^{\alpha}(K)\rightarrow \mathbf{V}(K)$ on the computational elements are
given by $\boldsymbol{\pi}_{K}=\mathbb{F}_{K}^{div}\circ\hat{\boldsymbol{\pi}}\circ[\mathbb{F}^{div}_{K}]^{-1}$.
Global projections $\boldsymbol{\Pi}_{h}:\mathbf{H}^{\alpha}(\Omega)\rightarrow\mathbf{V}_{h}$
are then piecewise defined: $(\boldsymbol{\Pi}_{h}\mathbf{q})|_{K}=\boldsymbol{\pi}_{K}(\mathbf{q}|_{K})$. Analogously, if $\hat{\lambda}$ denotes the $L^{2}$-projection on $\hat{U}$, a global projection $\Lambda: L^2(\Omega)\rightarrow U$  is piecewise defined:  $(\Lambda\varphi|_{K}={\lambda}_{K}(\varphi|_{K})$, where $\lambda_{K}=\mathbb{F}_{K}\circ\hat{\lambda}\circ\mathbb{F}_{K}^{-1}$.

The projections
$\hat{\boldsymbol{\pi}}$ and $\hat{\lambda}$ are required to verify the commutative De Rham property illustrated in the next  diagram 

\begin{equation}
\begin{array}[t]{ccc}
\mbox{\ensuremath{\mathbf{H}^{\alpha}(\hat{K})}} & \overset{\nabla\cdot}{\longrightarrow} & L^{2}(\hat{K})\\
\downarrow\hat{\boldsymbol{\pi}} &  & \downarrow\hat{\lambda}\\
\hat{\mathbf{V}} & \overset{\nabla\cdot}{\longrightarrow} & \hat{U}
\end{array}, \label{eq:deRham1}
\end{equation}
meaning that
\begin{equation}
\int_{\hat{K}}\nabla\cdot [\hat{\boldsymbol{\pi}}\mathbf{q}-\mathbf{q}]\varphi\;d\hat{K}=0,\;\;\forall\varphi\in \hat{U}.\label{eq:deRham2}
\end{equation}

\subsubsection*{Remarks}
 As suggested in \cite{DemkowiczHP} and applied in  \cite{Jay2005}, there is a general form to
define bounded projections $\hat{\boldsymbol{\pi}}$ of smooth functions $\mathbf{q}\in\mathbf{H}(div,\hat{K})$
onto $\hat{\mathbf{V}}$, verifying the local de Rham property (\ref{eq:deRham2}), without requiring any specific geometric aspect,  which is valid by  the usual vector approximation spaces, making use of their factorizations in terms of internal and edge components. They are defined by the requirements
\begin{eqnarray}
\int_{\partial \hat{K}}\hat{\boldsymbol{\pi}}\mathbf{q\cdot\boldsymbol{\eta}}\phi \;ds & = & \int_{\partial \hat{K}}\mathbf{q\cdot\boldsymbol{\eta}}\phi\; ds,\;\;\forall\phi\in P(\partial \hat{K}),\label{eq:Proj1}\\
\int_{\hat{K}}\nabla\cdot \hat{\boldsymbol{\pi}}\mathbf{q}\;\nabla \cdot \boldsymbol{\sigma}\;d\hat{K} & = & \int_{\hat{K}}\nabla \cdot \mathbf{q}\;\nabla\cdot\boldsymbol{\sigma}\;d\hat{K},\;\;\forall\boldsymbol{\sigma}\in\mathring{\hat{\mathbf{V}}},\label{eq:Proj2}\\
\int_{\hat{K}}\hat{\boldsymbol{\pi}}\mathbf{q} \cdot \boldsymbol{\sigma}\;d\hat{K} & = & \int_{\hat{K}}\mathbf{q} \cdot \boldsymbol{\sigma}\;d\hat{K},\;\;\forall\boldsymbol{\sigma}\in\mathring{\hat{\mathbf{V}}},\;\nabla \cdot \boldsymbol{\sigma}=0,\label{eq:Proj3}
\end{eqnarray}
where $P(\partial \hat{K})$ represents the space of normal traces of vector
functions in $\hat{\mathbf{V}}$. Note that equation
(\ref{eq:Proj2}) is trivially verified by divergence free functions
$\boldsymbol{\sigma}\in\mathring{\hat{\mathbf{V}}}$.
Therefore, it needs to be tested only for functions $\boldsymbol{\sigma}\in\mathring{\mathbf{V}}$
with $\nabla\cdot\boldsymbol{\sigma} =\varphi \neq 0$.
 
 As such, the projections $\hat{\boldsymbol{\pi}}$ admit a factorization $\hat{\boldsymbol{\pi}}  {\mathbf{q}} =  \hat{\boldsymbol{\pi}}^\partial {\mathbf{q}} + \mathring{\hat{\boldsymbol{\pi}}}{\mathbf{q}}$, in terms of edge and internal contributions, the  first  term $\hat{\boldsymbol{\pi}}^\partial {\mathbf{q}}$ being determined by the requirement (\ref{eq:Proj1}), the constraints (\ref{eq:Proj2})-(\ref{eq:Proj3}) determining the complementary component $\mathring{\hat{\boldsymbol{\pi}}} {\mathbf{q}}$.
The uniqueness of $\hat{\boldsymbol{\pi}}$ can be deduced from the conditions (\ref{eq:Proj1})-(\ref{eq:Proj3})
by assuming zero right hand sides 
\begin{eqnarray}
\int_{\partial \hat{K}}\hat{\boldsymbol{\pi}}\mathbf{q\cdot\boldsymbol{\eta}}\phi \;ds & = & 0,\;\;\forall\phi\in P(\partial \hat{K}),\label{eq:Proj1-0}\\
\int_{\hat{K}}\;\nabla\cdot \hat{\boldsymbol{\pi}}\mathbf{q}\;\nabla\cdot\boldsymbol{\sigma}\;d\hat{K} & = & 0,\;\;\forall\boldsymbol{\sigma}\in\mathring{\hat{\mathbf{V}}}(\hat{K}),\label{eq:Proj2-0}\\
\int_{\hat{K}}\hat{\boldsymbol{\pi}}\mathbf{q}\cdot\boldsymbol{\sigma}\;d\hat{K} & = & 0,\;\;\forall\boldsymbol{\sigma}\in\mathring{\hat{\mathbf{V}}}(K),\nabla\cdot\boldsymbol{\sigma}=0.\label{eq:Proj3-0}
\end{eqnarray}
Equation (\ref{eq:Proj1-0}) implies that $\hat{\boldsymbol{\pi}}\mathbf{q}\cdot\boldsymbol{\boldsymbol{\eta}}|_{\partial K}=0$,
meaning that $\hat{\boldsymbol{\pi}}\mathbf{q}\in\mathring{\hat{\mathbf{V}}}(K)$.
By taking $\boldsymbol{\sigma}=\hat{\boldsymbol{\pi}}\mathbf{q}$ in equations
(\ref{eq:Proj2-0}), it follows that $\nabla\cdot \hat{\boldsymbol{\pi}}\mathbf{q}=0$.
Therefore, $\boldsymbol{\sigma}=\hat{\boldsymbol{\pi}}\mathbf{q}$ can be applied
in equation (\ref{eq:Proj3-0}) to conclude that $\hat{\boldsymbol{\pi}}\mathbf{q}=0$.

In order to verify the  commutative relation (\ref{eq:deRham2}), 
 first note that it is valid for constant $\varphi=1$ as a consequence of  (\ref{eq:Proj1}),
\begin{eqnarray*}
\int_{\hat{K}}\nabla\cdot [\hat{\boldsymbol{\pi}}\mathbf{q}-\mathbf{q}]\;d\hat{K}&=&\int_{\partial \hat{K}}[\hat{\boldsymbol{\pi}} \mathbf{q}-\mathbf{q}]\cdot\hat{\boldsymbol{\eta}}\; ds \\
&=&\int_{\partial \hat{K}}[\hat{\boldsymbol{\pi}}^\partial {\mathbf{q}}-\mathbf{q}]\cdot\hat{\boldsymbol{\eta}}\; ds =0.\\
\end{eqnarray*}
Therefore, it remains to prove (\ref{eq:deRham2}) for $\varphi \in \hat{U}$
with zero mean, $\varphi \neq 0$.  Let $\boldsymbol{\sigma}\in\mathring{\hat{\mathbf{V}}}$ be such
that $\nabla\cdot\boldsymbol{\sigma}=\varphi$. Thus, by (\ref{eq:Proj2}), the desired result holds
\begin{eqnarray*}
\int_{\hat{K}}\nabla\cdot [\hat{\boldsymbol{\pi}}\mathbf{q}-\mathbf{q}] \varphi \;d\hat{K}&=&\int_{\hat{K}}\nabla\cdot [\hat{\boldsymbol{\pi}}\mathbf{q}-\mathbf{q}]\nabla \cdot \boldsymbol{\sigma} \;d\hat{K} =0.
\end{eqnarray*}

\clearpage

\subsection*{Projection errors}
The following projection error estimates hold for spaces based on shape-regular 
meshes  ${\cal T}_{h}$,  and  for sufficiently smooth vector functions $\mathbf{q}$ and scalar functions $u$  (see Theorem 4.1 and Theorem 4.2 in \cite{ABF}, and  Theorem 1 in \cite{ABF2002})
\begin{eqnarray}
||\mathbf{q}-\boldsymbol{\Pi}_{h}\mathbf{q}||_{\mathbf{L}^{2}(\Omega)}&\lesssim & \;h^{k+1}||\mathbf{q}||_{\mathbf{H}^{k+1}(\Omega)},\label{eq:projflux}\\
||\nabla\cdot(\mathbf{q}-\boldsymbol{\Pi}_{h}\mathbf{q})||_{L^{2}(\Omega)}&\lesssim & \;h^{r+1}||\nabla\cdot\mathbf{q}||_{H^{r+1}(\Omega)}, \label{eq:errorprojdiv}\\
||u-\Lambda_{h}u||_{L^{2}(\Omega)}&\lesssim & \;h^{t+1}||u||_{H^{t+1}},\label{eq:errorProjPot}
\end{eqnarray}
where  $k$  and $r$ are determined by the total degree of polynomials such that $ [\mathbb{P}_k]^2 \subset \mathbf{V}(K)$, and  $ \mathbb{P}_r \subset \nabla \cdot \mathbf{V}(K)$, and $t$ is such that $ \mathbb{P}_{t} \subset U(K)$. Here and in what follows we write $a\lesssim b$ whenever $a \leq C \,b$ for a constant $C$ not depending on essential quantities. The leading constants on (\ref{eq:projflux}) and (\ref{eq:errorprojdiv}) depends only on the bound for  corresponding projection $\hat{\boldsymbol{\pi}}$ on the master element,  and on the shape regularity constant of ${\cal T}_{h}$.\\


\subsubsection*{Optimal conditions}

The  parameters $k$ and $r$ determining the convergence rates in (\ref{eq:projflux}) and (\ref{eq:errorprojdiv}) can be easily obtained   when affine transformations $F_K$ are used, since for them the action of the Piola transformation  preserves the polynomial vector fields  of the reference space $\hat{\mathbf{V}}$. For instance, in the case of $BDM_k$ spaces for triangles, 
$\hat{\mathbf{V}}_{BDM_k}=[\mathbb{P}_k]^2$,  and $\nabla\cdot \hat{\mathbf{V}}_{BDM_k}=\mathbb{P}_{k-1}$.  Similarly, for Raviart-Thomas spaces  on affine quadrilaterals, $\hat{\mathbf{V}}_{RT_k}=\mathbb{Q}_{k+1,k}\times \mathbb{Q}_{k,k+1}\supset [\mathbb{P}_{k}]^2$,  and $\nabla\cdot \hat{\mathbf{V}}_{RT_k}=\mathbb{Q}_{k,k}\supset \mathbb{P}_{k}$.
However, for general  quadrilateral elements,  with non-constant Jacobian determinants, this is a more subtle task,  as discussed in \cite{ABF}. For that,  the following optimal concepts are helpful, depending on the  geometric mapping $F_K$:\\

\begin{itemize}
\item Definition A: $\hat{\mathbf{S}}_{k}^{K},\;k\geq 1$, is the vector polynomial space on $\hat{K}$ of minimal dimension such that the following property holds 
\begin{equation}
\hat{\mathbf{S}}_{k}^K \subset \hat{\mathbf{V}}  \Longleftrightarrow \mathbf{V}(K)=\mathbb{F}_{K}^{div}\hat{\mathbf{V}} \supset[\mathbb{P}_{k}]^{2}. \label{eq:opt1}
\end{equation}

\item Definition B: $\hat{R}_{r}^K,\;r\geq 1$, is the scalar polynomial space on $\hat{K}$ of minimal
dimension  such that the following property holds 
\begin{equation}
\hat{{R}}_{r}^K \subset \nabla\cdot\hat{\mathbf{V}}  \Longleftrightarrow \nabla\cdot\mathbf{V}(K)=\nabla\cdot\mathbb{F}_{K}^{div}\hat{\mathbf{V}} \supset\mathbb{P}_{r}. \label{eq:opt3}
\end{equation}
\end{itemize}

The characterization of the optimal spaces  $\hat{\mathbf{S}}_{k}^K$ and $\hat{R}_{r}^K$ has been established in \cite{ABF} for bilinearly mapped quadrilateral elements $K$.  Precisely, 
\begin{itemize}
\item $\hat{\mathbf{S}}_{k}^{K},\;k\geq 1$, is the subspace of codimension one in $\hat{\mathbf{V}}_{RT_k}=\mathbb{Q}_{k+1,k}\times \mathbb{Q}_{k,k+1}$, spanned by its vector fields, except $(\hat{x}^{k+1}\hat{y}^k,0)$ and $(0, \hat{x}^{k}\hat{y}^{k+1})$, which are replaced by the single vector $(\hat{x}^{k+1}\hat{y}^k, -\hat{x}^{k}\hat{y}^{k+1})$ (see \cite{ABF}, Theorem 3.1).
\item  $\hat{R}_{r}^K,\;r\geq 1$, is the subspace of
codimension one of $\mathbb{Q}_{r+1,r+1}$,  spanned by all the monomials, except $\hat{x}^{r+1}\hat{y}^{r+1}$
(see \cite{ABF}, Theorem 3.2).
\end{itemize}
The $ABF_k$ spaces are then defined as  enrichment of the  $RT_k$ spaces, by taking $\hat{\mathbf{V}}_{ABF_k}=\mathbb{Q}_{k+2,k}\times \mathbb{Q}_{k,k+2}\supset \hat{\mathbf{V}}_{RT_k} \supset \hat{\mathbf{S}}_{k}^{K}$, whose divergence space is $\nabla\cdot \hat{\mathbf{V}}_{ABF_k}= \hat{R}_{r}^K$, ensuring convergence rates of order $k+1$ for flux and flux divergence for approximations based on general bilinearly mapped quadrilaterals, enhancing the divergence accuracy of $RT_k$ approximations  based on such non-affine meshes (of order $k$). 


\section{Enriched  approximation space configurations}\label{section3}

Suppose that vector polynomial spaces $\hat{\mathbf{V}}_k$ are defined on the master element $\hat{K}$, for which  the following properties hold:
\begin{enumerate}
\item A direct factorization $\hat{\mathbf{V}}_k=\hat{\mathbf{V}}^\partial_k  \oplus \mathring{\hat{\mathbf{V}}}_k$ is defined in terms  of  edge and internal flux functions, the edge functions having normal components  of degree $k$ over $\partial \hat{K}$.
\item $\hat{\mathbf{V}}_k$ contains the optimal vector polynomial space $\hat{\mathbf{S}}_k^K$.
\item The associated divergence space $\hat{U}_k = \nabla \cdot \hat{\mathbf{V}}_k$  contains the optimal scalar polynomial space $\hat{R}_r^K$, and $\mathbb{F}_K(\hat{U}_k)$ contains $\mathbb{P}_t$.
\item Bounded projections
$\hat{\boldsymbol{\pi}}_k: \mathbf{H}^{\alpha}(\hat{K})\rightarrow \hat{\mathbf{V}}_k$ and $\hat{\lambda}_k: L^2(\hat{K})\rightarrow \hat{U}_k$  verifying the commutative De Rham property are available.
\item  The projections
$\hat{\boldsymbol{\pi}}_k$ can be factorized as $\hat{\boldsymbol{\pi}}_k  \hat{\mathbf{q}} = \hat{\boldsymbol{\pi}}_k^\partial \hat{\mathbf{q}} + \mathring{\hat{\boldsymbol{\pi}}}_k\hat{\mathbf{q}}$, in terms of edge and internal contributions, which are determined by constraints, as described in (\ref{eq:Proj1})-(\ref{eq:Proj3}).
\end{enumerate}

Under these conditions,  the desired enriched versions $\hat{\mathbf{V}}^{n+}_k$, $n\geq 1$, are defined  by adding  to  $\hat{\mathbf{V}}_k$ higher degree internal shape functions of the original space at level $k+n$,  while keeping the original border fluxes at level $k$.
\[\hat{\mathbf{V}}_{k}^{n+} = \hat{\mathbf{V}}^{\partial}_{k}\oplus\mathring{\hat{\mathbf{V}}}_{k+n}.\]
As the edge components $\hat{\mathbf{V}}^{\partial}_{k}$ are kept fixed, the optimal space contained in $\hat{\mathbf{V}}_{k}^{n+}$ is not improved, remaining as $\hat{\mathbf{S}}_k^K$. But, the corresponding enriched divergence spaces are now
 \[\hat{U}_k^{n+} = \nabla \cdot \hat{\mathbf{V}}_{k}^{n+} = \hat{U}_{k+n},\]
 containing  the optimal  space $\hat{R}_{r+n}^K$, and such that $ \mathbb{P}_{t+n} \subset \mathbb{F}_K(\hat{U}_k^{n+})$,  a  guarantee for higher divergence accuracy  and enhanced potential approximations for schemes based on these frameworks.


\subsection{Projections}
Projections $\hat{\boldsymbol{\pi}}^{n+}_k:\mathbf{H}^{\alpha}(\hat{K})\rightarrow\hat{\mathbf{V}}_k^{n+}$ are defined as $\hat{\boldsymbol{\pi}}^{n+}_k{\mathbf{q}} = \hat{\boldsymbol{\pi}}_k^{n+,\partial} \hat{\mathbf{q}} + \mathring{\hat{\boldsymbol{\pi}}}_k^{n+}\hat{\mathbf{q}}$, where 
\begin{enumerate}
\item  The edge component $\hat{\boldsymbol{\pi}}_k^{n+,\partial} \hat{\mathbf{q}}=\hat{\mathbf{q}}^\partial \in \hat{\mathbf{V}}^{\partial}$ is determined
by 
\begin{equation}\int_{\partial\hat{K}}[{\mathbf{q}}- \hat{\mathbf{q}}^\partial ]\cdot\hat{\boldsymbol{\eta}}\;\phi ds=0,\;\forall\phi\in P_{k}(\partial\hat{K}), \label{eq:Proj1+}
\end{equation}
$P_{k}(\partial\hat{K})$ representing the normal traces of functions in $\hat{\mathbf{V}}_k$.

\item The internal term $\mathring{\hat{\boldsymbol{\pi}}}_k^{n+}\hat{\mathbf{q}}=\mathring{\hat{\mathbf{q}}} \in\mathring{\hat{\mathbf{V}}}_k^{n+}=\mathring{\hat{\mathbf{V}}}_{k+n}$ is taken
as $\mathring{\hat{\mathbf{q}}} =\mathring{\hat{\boldsymbol{\pi}}}_{k+n}({\mathbf{q}}-\mathbf{q}_{\partial})$,
where $\mathring{\hat{\boldsymbol{\pi}}}_{k+n}$ is the internal projection component of the original scheme  at level $k+n$.  Precisely,  it is determined by the following constraints, valid for $\boldsymbol{\sigma}\in\mathring{\hat{\mathbf{V}}}_{k+n}$,
\begin{align}
\int_{\hat{K}}\nabla \cdot \mathring{\hat{\mathbf{q}}}\;\nabla \cdot \boldsymbol{\sigma}d\hat{K} & =  \int_{\hat{K}}\nabla \cdot (\mathbf{q} - \hat{\mathbf{q}}^\partial)\;\nabla\cdot\boldsymbol{\sigma}\;d\hat{K},\label{eq:Proj2+}\\
\int_{\hat{K}}\mathring{\hat{\mathbf{q}}}\cdot \boldsymbol{\sigma}d\hat{K} & =  \int_{\hat{K}}(\mathbf{q}-\hat{\mathbf{q}}^\partial) \cdot \boldsymbol{\sigma}d\hat{K}, \mbox{with} \;\nabla \cdot \boldsymbol{\sigma}=0.\label{eq:Proj3+}
\end{align}
\end{enumerate}

For the verification of  uniqueness and the the corresponding  De Rham commutative relation 
\begin{equation}
\int_{\hat{K}}\nabla\cdot [\hat{\boldsymbol{\pi}}_k^{n+}\mathbf{q}-\mathbf{q}]\varphi\;d\hat{K}=0,\;\;\forall\varphi\in \hat{U}_{k+n},\label{eq:deRham2+}
\end{equation}
the same arguments can be used, as done for the  original scheme.

\subsection*{Projection errors}

Recall that the basic hypotheses on the original space  imply that,  on the mapped elements $K\in \mathcal{T}_h$,  the enriched scalar  spaces verify $\mathbb{P}_{t+n} \subset \mathbb{F}_K(\hat{U}_k^{n+})$ and,  the local vector spaces are such that
\begin{eqnarray*}
[\mathbb{P}_{k}]^{2} & \subset &\mathbf{V}^{n+}_k(K)  \\
  \mathbb{P}_{r+n} & \subset & \nabla\cdot\mathbf{V}_k^{n+}(K).
\end{eqnarray*}
Therefore, for the  globally extended projections  $\boldsymbol{\Pi}_{h,k}^{n+}:\mathbf{H}^{\alpha}(\Omega)\rightarrow\mathbf{V}_{h,k}^{n+}$,  and $\Lambda_{h,k}^{n+}:L^2\Omega)\rightarrow U_{h,k}^{n+}=U_{h,k+n}$,  the  error estimations  (\ref{eq:projflux})-(\ref{eq:errorProjPot}) are translated into
\begin{eqnarray}
||\mathbf{q}-\boldsymbol{\Pi}_{h,k}^{n+}\mathbf{q}||_{\mathbf{L}^{2}(\Omega)}&\lesssim & \;h^{k+1}||\mathbf{q}||_{\mathbf{H}^{k+1}(\Omega)},\label{eq:projfluxnplus}\\
||\nabla\cdot(\mathbf{q}-\boldsymbol{\Pi}_{h,k}^{n+}\mathbf{q})||_{L^{2}(\Omega)}&\lesssim & \;h^{r+n+1}||\nabla\cdot\mathbf{q}||_{H^{r+n+1}(\Omega)}, \label{eq:errorprojdivnplus}\\
||u-\Lambda_{h,k}^{n+}u||_{L^{2}(\Omega)}&\lesssim & \;h^{t+n+1}||u||_{H^{t+n+1}}.\label{eq:errorProjPotnplus}
\end{eqnarray}

\section{Application to mixed finite element formulations}\label{section4}
Consider the  problem  of finding $\boldsymbol{\sigma}\in \mathbf{H}(\mbox{div},\Omega)$, and $u\in L^2(\Omega)$ satisfying $\nabla \cdot\boldsymbol{\sigma} = f$,  $\boldsymbol{\sigma}  =-\mathcal{K}\nabla u$, and
$u|_{\partial \Omega}  = u_{D}$, where  $f\in L^2(\Omega)$, and $u_D\in  H^{1/2}(\partial \Omega)$. The tensor $\mathcal{K}$ is assumed to be a symmetric positive-definite matrix, composed by functions in $L^\infty(\Omega)$. Given approximation spaces $\mathbf{V}_{h}\subset \mathbf{H}(\mbox{div},\Omega)$ for the variable $\boldsymbol{\sigma}$,  and $U_h\subset L^2(\Omega)$ for the variable $u$,  based on  the family of meshes ${\cal T}_{h}$ of $\Omega$,    consider the  discrete variational mixed formulation for this problem \cite{BrezziFortin1991}: Find $(\boldsymbol{\sigma}_{h},u_{h})\in\left(\mathbf{V}_{h}\times U_{h}\right)$, such that $\forall\mathbf{q}\in\mathbf{V}_{h}$,
and $\forall\varphi\in U_{h}$
\begin{align}
\int_\Omega \mathcal{K}^{-1}\boldsymbol{\sigma}_h \cdot\mathbf{q}\; d\Omega - \int_\Omega u_h \nabla\cdot \mathbf{q} \;d\Omega & = -<u_{D},\mathbf{q}\cdot\boldsymbol{\eta}>, \label{eq:MF1-1}\\
\int_\Omega \varphi \nabla\cdot \boldsymbol{\sigma}_{h} \;d\Omega &= \int_\Omega f \varphi \;d\Omega,\label{MF2-1}
\end{align}
where $\boldsymbol{\eta}$ is the outward unit normal vector on $\partial \Omega$, and $<\cdot ,\cdot>$ represents the duality pairing of $H^{1/2}(\partial \Omega)$ and $H^{-1/2}(\partial \Omega)$.

\subsection*{A priori error estimations}
Consider an enriched space configuration  $\mathbf{V}_{h}= \mathbf{V}_{h,k}^{n+}$ and $U_h=U_{h,k}^{n+}$, $n\geq 1$, based on shape regular partitions ${\cal T}_{h}$  of a convex region  $\Omega$, as described in the previous sections.  Using  the general estimations in \cite{ABF}, Theorem 6.1, in terms of projection errors of the exact solutions, and assuming  sufficiently regular exact solutions,  the following error estimations hold 
\begin{align}
||\boldsymbol{\sigma}-\boldsymbol{\sigma}_{h}||_{\mathbf{L}^{2}(\Omega)} & \lesssim h^{k+1}||\boldsymbol{\sigma}||_{\mathbf{H}^{s+1}(\Omega)},\label{eq:errorfluxmixed}\\
||\nabla\cdot(\boldsymbol{\sigma}-\boldsymbol{\sigma}_{h})||_{L^{2}(\Omega)} & \lesssim h^{r+n+1}\;||\nabla\cdot\boldsymbol{\sigma}||_{H^{r+n+1}},\label{eq:errordivmixed}\\
||u-u_{h}||_{L^{2}(\Omega)} & \lesssim h^{q+1}||u||_{H^{q+1}},\label{eq:errorpotentialmixed}
\end{align}
where $q=\min \{k+1,r+n+1,t+n\}$.  The   estimations for the flux  in (\ref{eq:errorfluxmixed})
and for the flux divergence in (\ref{eq:errordivmixed}) are derived directly from the  projection errors (\ref{eq:projflux})  and (\ref{eq:errorprojdiv}). For the potential variable $u$, the order of accuracy in (\ref{eq:errorpotentialmixed}) follows from  similar arguments as in the proof of Theorem 6.2 in \cite{ABF}, and from the projection error  (\ref{eq:errorProjPot}). The convexity of $\Omega$ only plays a role for the elliptic regularity property, used to get (\ref{eq:errorpotentialmixed}). Note that, as $n$ increases, the best potential  accuracy is  $k+2$, one unit more than the flux accuracy. But,   divergence  accuracy can be enhanced as desired, with  rates of order $r+n+1$.

Table \ref{tab:orders} illustrates  the convergence orders of the errors in $L^2$-norms  for  flux,  potential, and  flux divergence variables obtained when the enriched    spaces configurations  $\mathbf{V}_{h,k}^{n+} \; U_{h,k}^{n+}$ are  used for the discretization of mixed finite element formulation for the Darcy's model problem, corresponding to the particular original spaces  $BDM_{k}$ based on triangles ( $t=r=k-1$, and $s=k$), and  $RT_k$  for quadrilaterals  ( $s=t=k$, and $r=k$ for affine meshes, or $r=k-1$ otherwise). Recall that the enriched space $BDM_{k}^{+}$ corresponds to the $BDFM_{k+1}$ case.

\begin{table}[htb]
\begin{centering}
\begin{tabular}{|c|c|c|c|c|c|c|c|}
\hline 
\multicolumn{2}{|c|}{} & \multicolumn{2}{c|}{{\footnotesize{}Flux}} & \multicolumn{2}{c|}{{\footnotesize{}Potential}} & \multicolumn{2}{c|}{{\footnotesize{}Divergence}}\tabularnewline \hline
 {\footnotesize{}Element} &  {\footnotesize{}Space} &{\footnotesize{}A} & {\footnotesize{}N-A} & {\footnotesize{}A} & {\footnotesize{}N-A} & {\footnotesize{}A} & {\footnotesize{}N-A}\tabularnewline
\hline 
\multirow{3}{*}{{\footnotesize{}${\cal T}$}} & {\scriptsize{}$BDM_{k}$} & {\scriptsize{}$k+1$} & {\scriptsize{}-} & {\scriptsize{}$k$} & {\scriptsize{}-} & {\scriptsize{}$k$} & {\scriptsize{}-}\tabularnewline
& {\scriptsize{}$BDFM_{k+1}$} & {\scriptsize{}$k+1$} & {\scriptsize{}-} & {\scriptsize{}$k+1$} & {\scriptsize{}-} & {\scriptsize{}$k+1$} & {\scriptsize{}-}\tabularnewline
 & {\scriptsize{}$\mathbf{V}_{h,k}^{n+},\; U_{h,k}^{n+}$, $n\geq2$} & {\scriptsize{}$k+1$} & {\scriptsize{}-} & {\scriptsize{}{$k+2$}} & {\scriptsize{}-} & {\scriptsize{}{$k+n+1$}} & {\scriptsize{}-}\tabularnewline
\hline 
\multirow{2}{*}{{\footnotesize{}${\cal R}$}} & {\scriptsize{}$RT_k$} & {\scriptsize{}$k+1$} & {\scriptsize{}$k+1$} & {\scriptsize{}$k+1$} & {\scriptsize{}$k+1$} & {\scriptsize{}$k+1$} & {\scriptsize{}{$k$}}\tabularnewline
 & {\scriptsize{}$\mathbf{V}_{h,k}^{n+},\; U_{h,k}^{n+}$, $n\geq 1$} & {\scriptsize{}$k+1$} & {\scriptsize{}$k+1$} & {\scriptsize{}{$k+2$}} & {\scriptsize{}{$k+2$}} & {\scriptsize{}$k+n+1$} & {\scriptsize{}{$k+n$}}\tabularnewline
\hline 
\end{tabular}
\par\end{centering}
\caption{\label{tab:orders} Convergence orders in $L^2$-norms of approximate solutions for the mixed formulation using  space configurations $BDM_{k}$ for triangles, and $RT_k$ for quadrilaterals, as original spaces, and  their enriched versions $\mathbf{V}_{k,h}^{n+}\; U_{k,h}^{n+}$, $n\geq 1$, using affine (A) or non-affine (N-A) meshes $\mathcal{T}_h$.} 
\end{table}

\section{Numerical results}\label{section5}
In this section we  illustrate the approximation results for the enriched versions of type $\mathbf{V}_{k,h}^{n+} \; U_{h,k}^{n+}$ proposed in the previous sections.  By taking as original spaces $RT_k$ based on  quadrilaterals, and $BDM_{k}$ for triangles,   the simplified notation $RT_k^{n+}$ and $BDM_{k}^{n+}$ is adopted for  identifying  their  enriched versions.  These space configurations are used for the discretization of the mixed finite element formulation for a Darcy's model problem defined in $\Omega=(0,1)\times(0,1)$,  the load function 
$f$ being chosen such that the   exact  solution is  
$u=\frac{\pi}{2} - \arctan\left[5 \left(\sqrt{(x - 1.25)^2 +(y+0.25)^2}-\frac{\pi}{3}\right)\right]$.

The partitions ${\cal T}_{h}$ of $\Omega$ are the same ones used  in tests shown in \cite{Hdiv2D2017}. Namely, uniform rectangular meshes are considered with spacing $h = 2^{-i},$ $i=2,\cdots,5$, and   triangular meshes are
constructed from them by diagonal subdivision. In  the trapezoidal meshes, the elements have a basis of length $h$, and vertical parallel sides   of lengths $0.75h$ and $1.25h$.  

\begin{figure}
\begin{centering}
\begin{tabular}{cc}
Rectangular meshes & Triangular meshes \\
\includegraphics[scale=0.35]{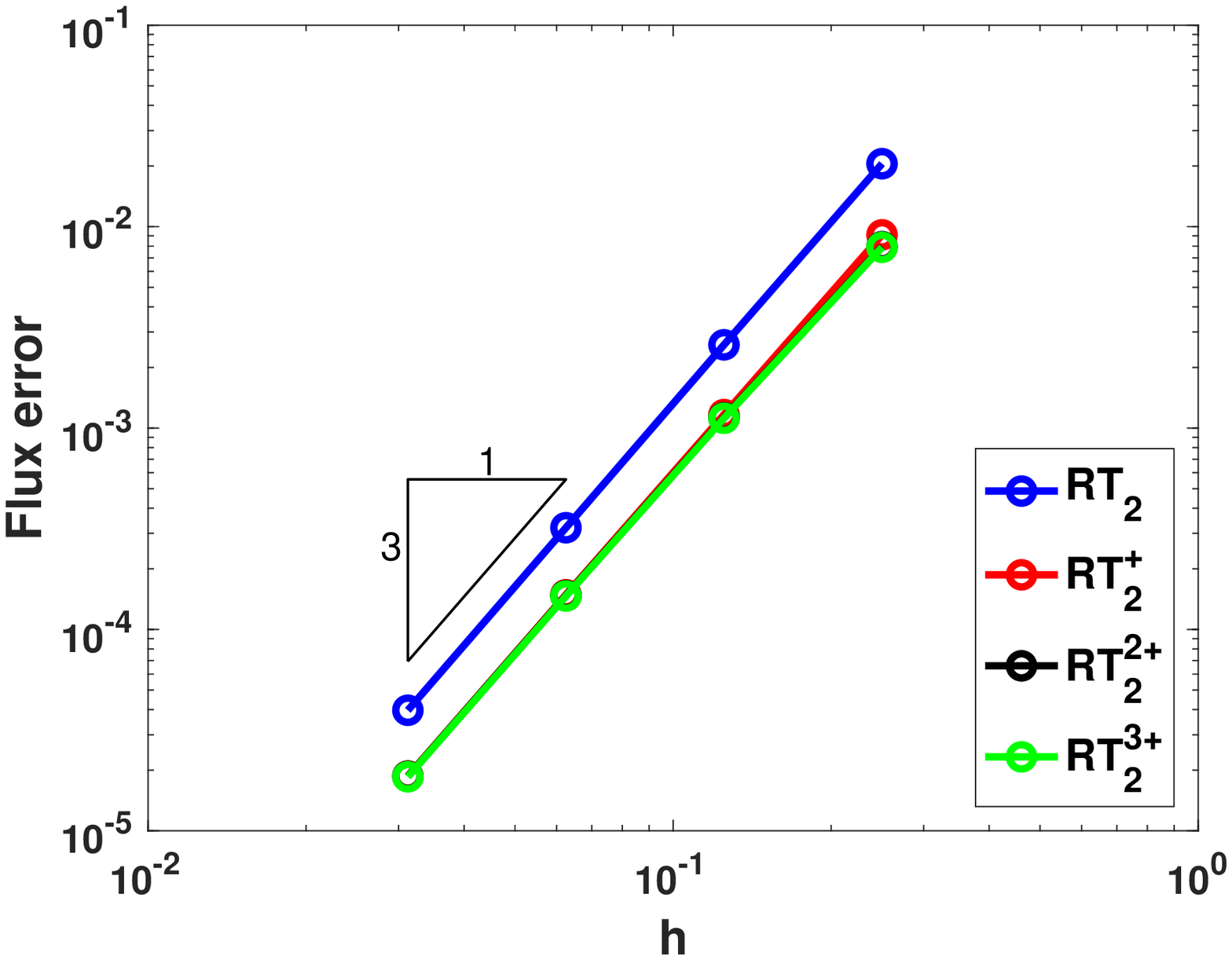} & \includegraphics[scale=0.35]{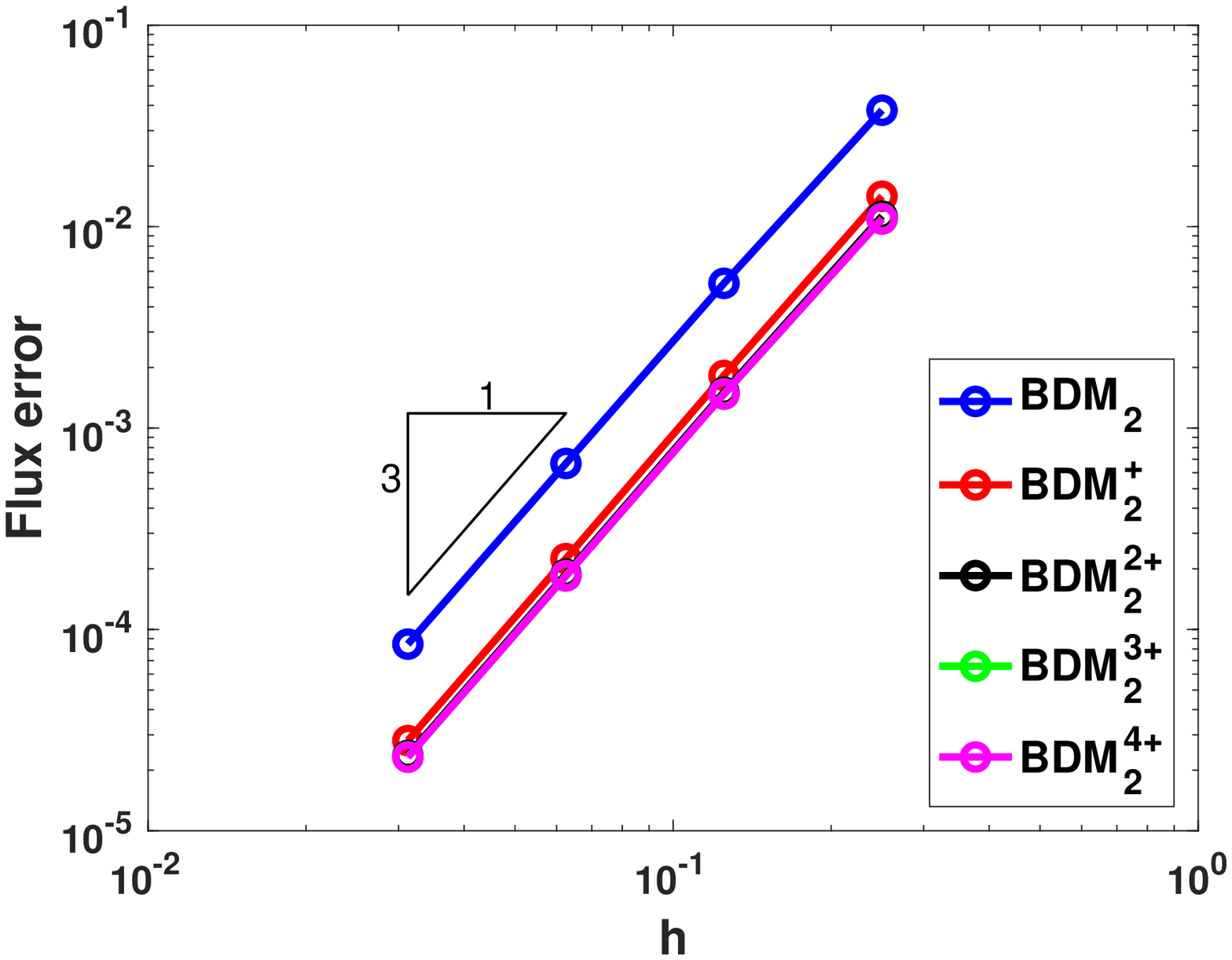}\tabularnewline
\includegraphics[scale=0.35]{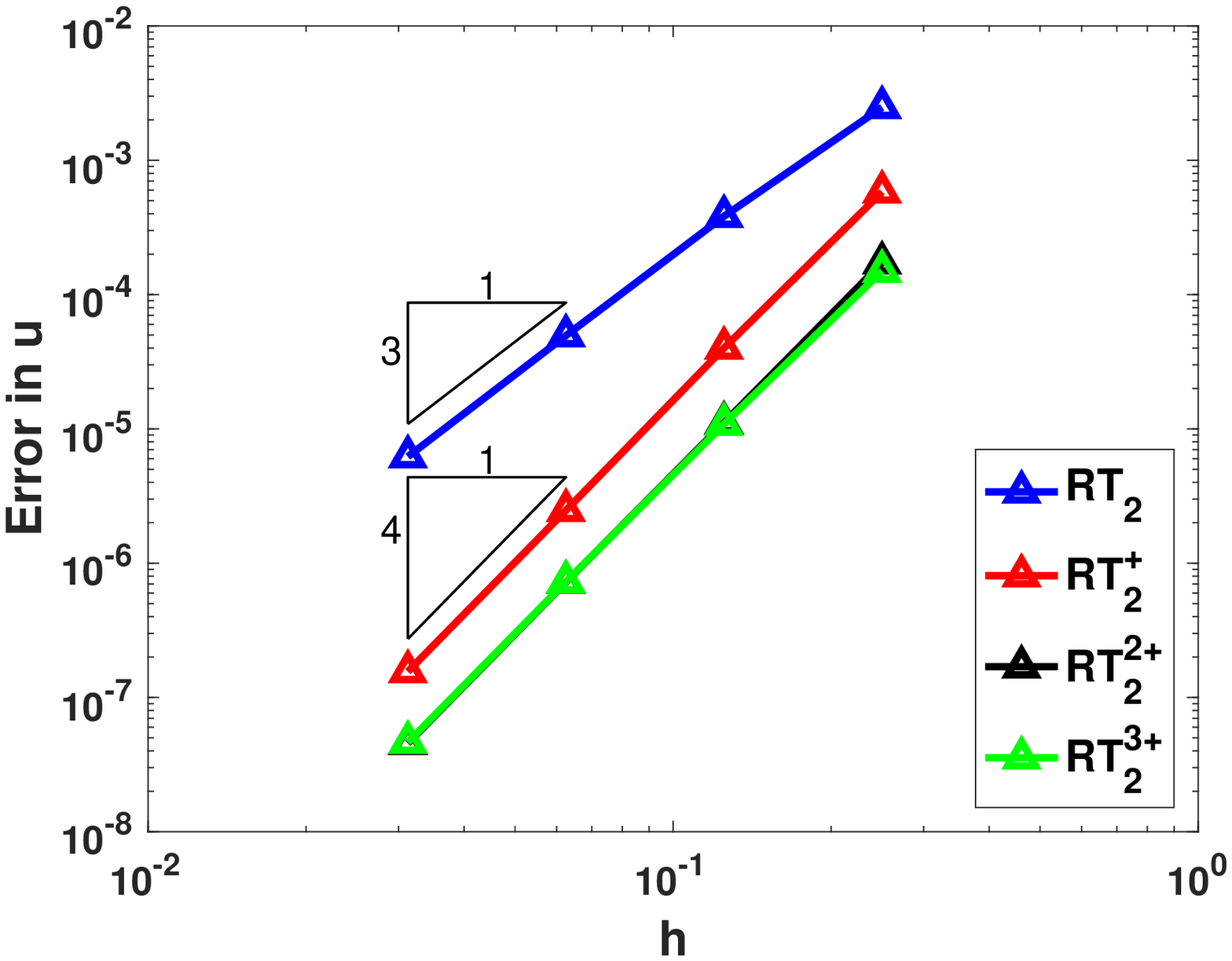} & \includegraphics[scale=0.35]{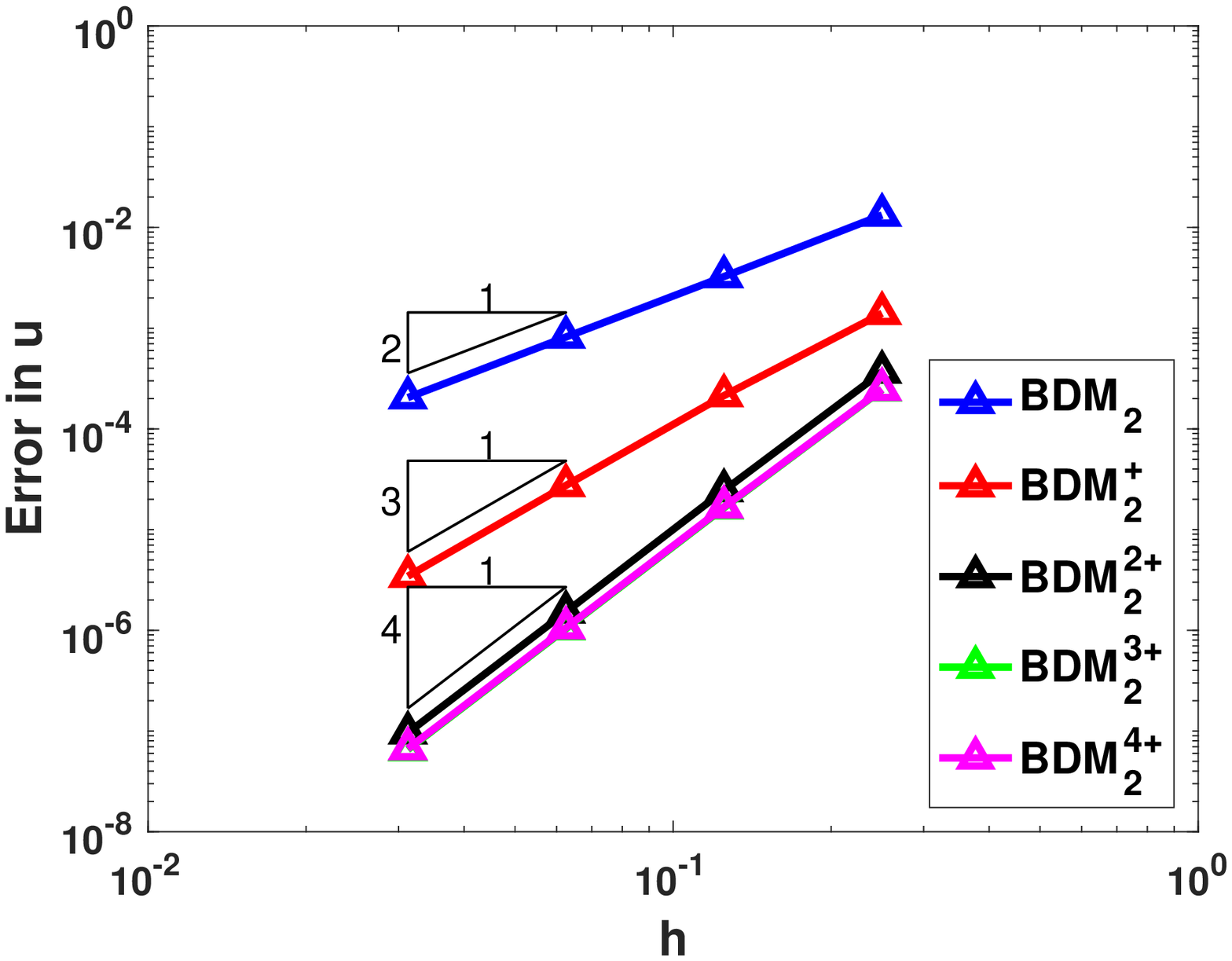}\tabularnewline
\includegraphics[scale=0.35]{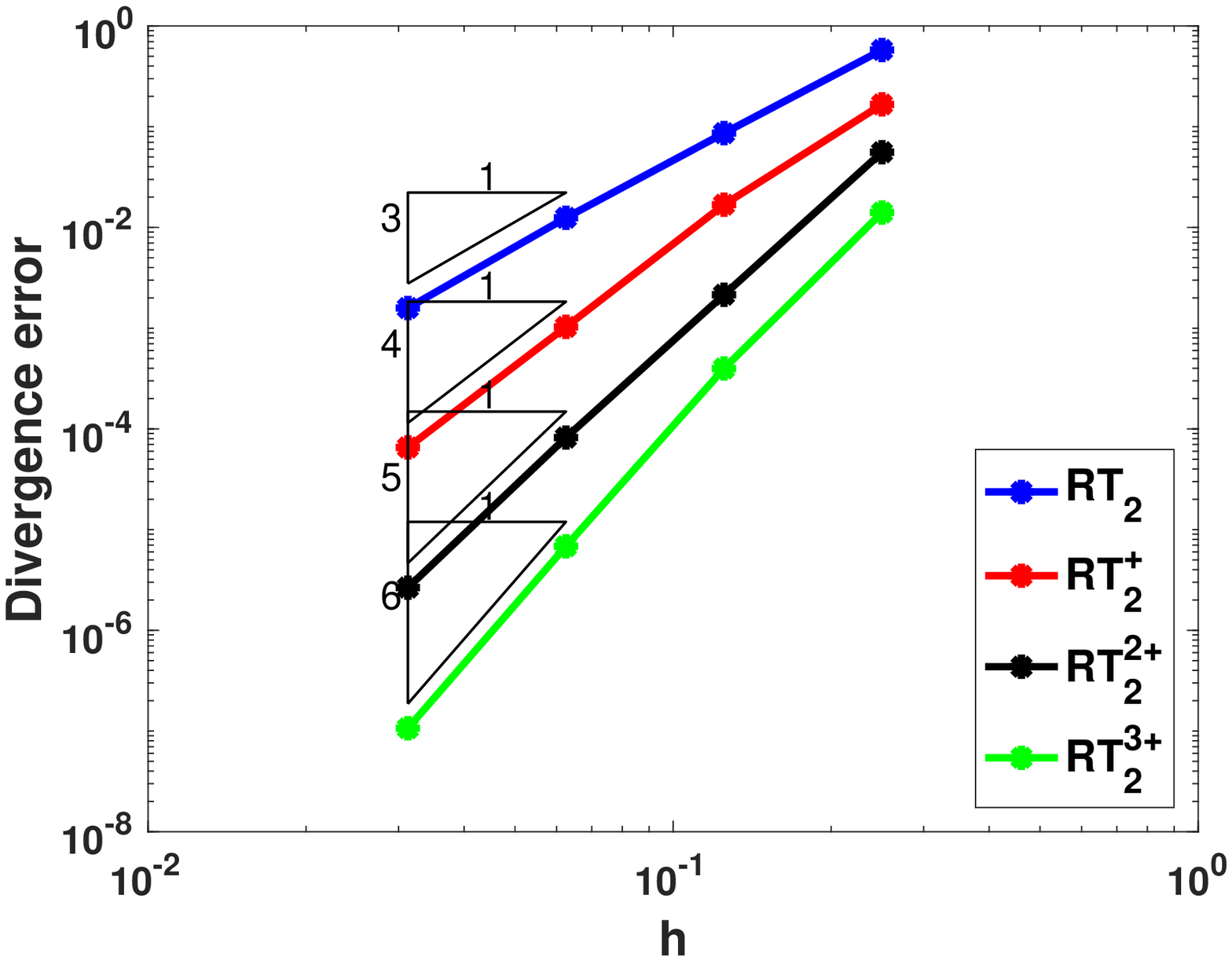} & \includegraphics[scale=0.35]{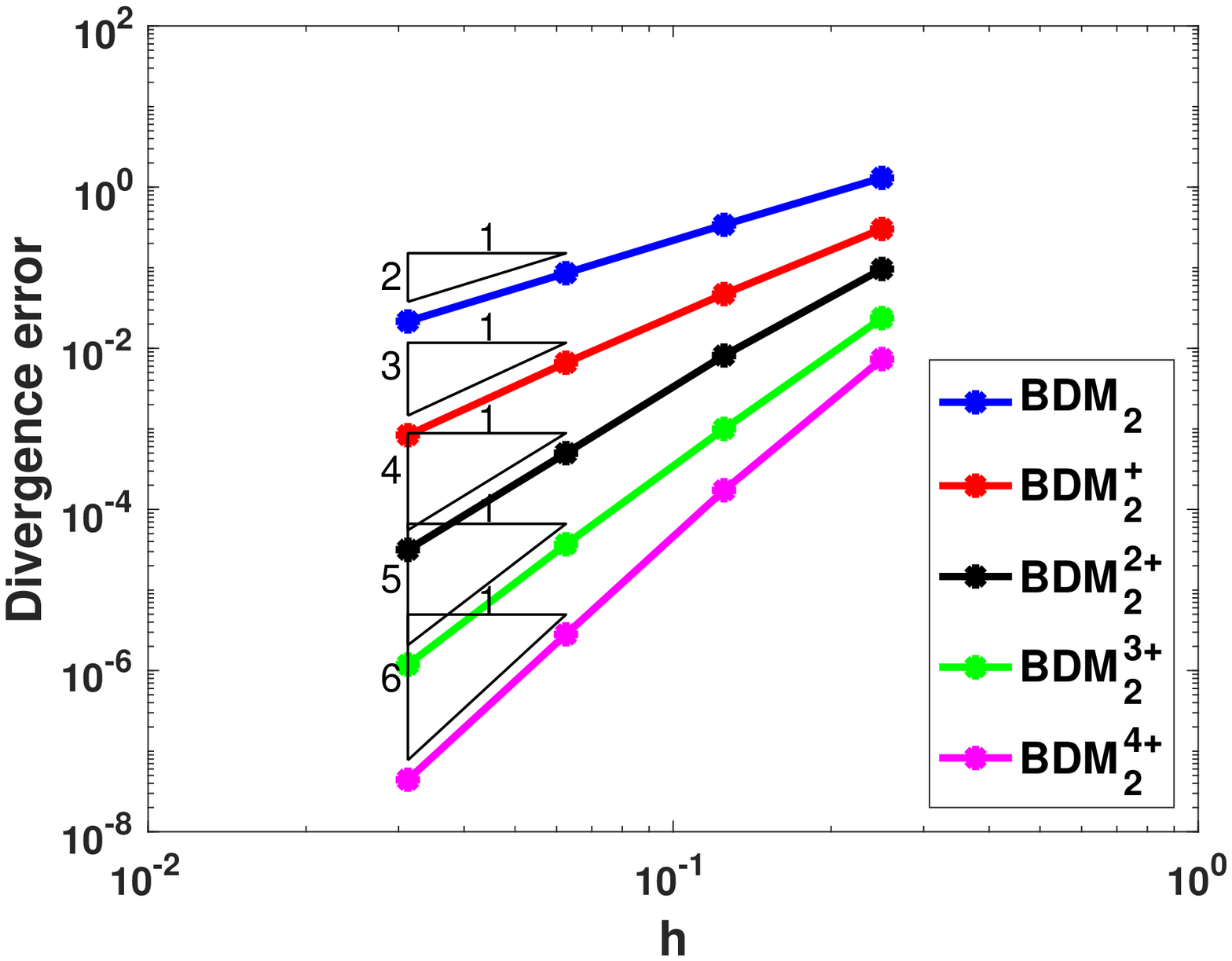}\tabularnewline
\end{tabular}
\par\end{centering}
\caption{$L^2$-errors in the flux (top), in $u$ (middle) and divergence (bottom), using approximation space configurations $RT_2^{n+}$ based rectangular meshes (right)  and  
$BDM_{2}^{n+}$ for triangles (left).\label{fig:Errors-QTAffinek2}}
\end{figure}

Figure \ref{fig:Errors-QTAffinek2} illustrates the cases of  approximation space configurations $RT_2^{n+}$ based rectangular meshes   and  
$BDM_{2}^{n+}$ for  triangular ones. The expected rates of convergence shown in Table \ref{tab:orders} are verified.  It should also be observed that stronger space enrichment  has practically no effect on the error values of the flux, for $n\geq 1$. Similarly,   the error of the potential remains almost the same  after applying the second enrichment. 

\begin{figure}
\begin{centering}
\begin{tabular}{cc}
\includegraphics[scale=0.35]{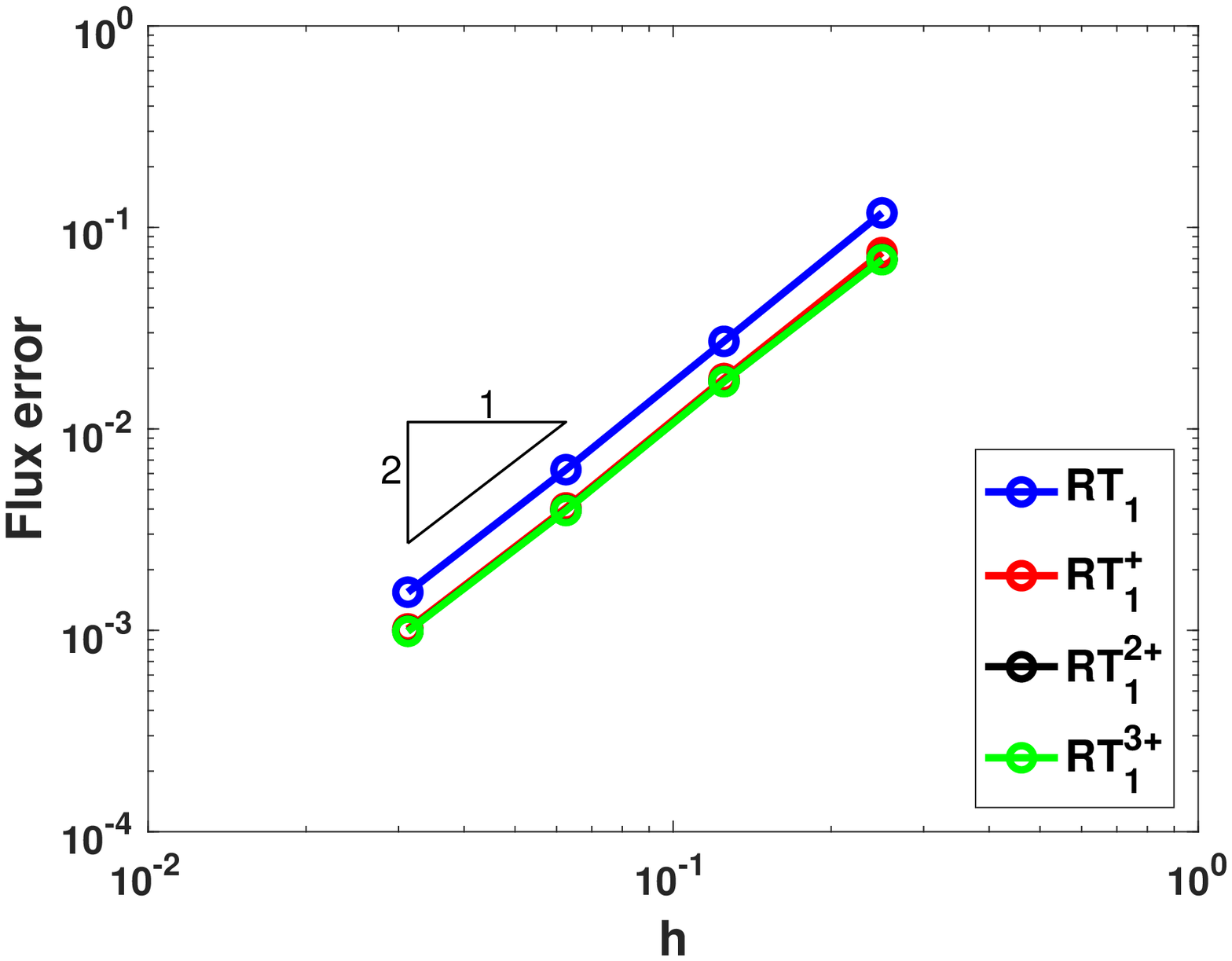} & \includegraphics[scale=0.35]{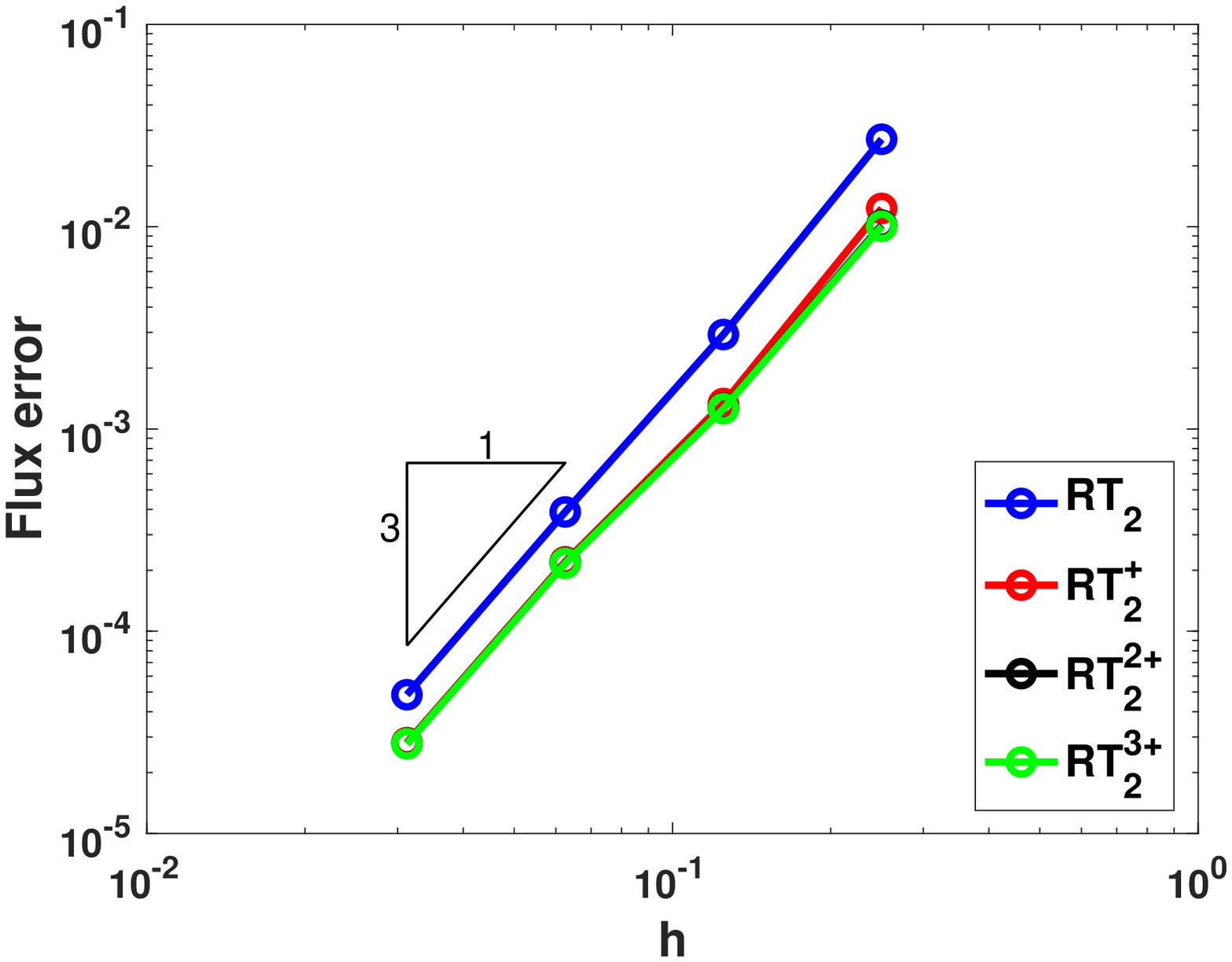}\tabularnewline
\includegraphics[scale=0.35]{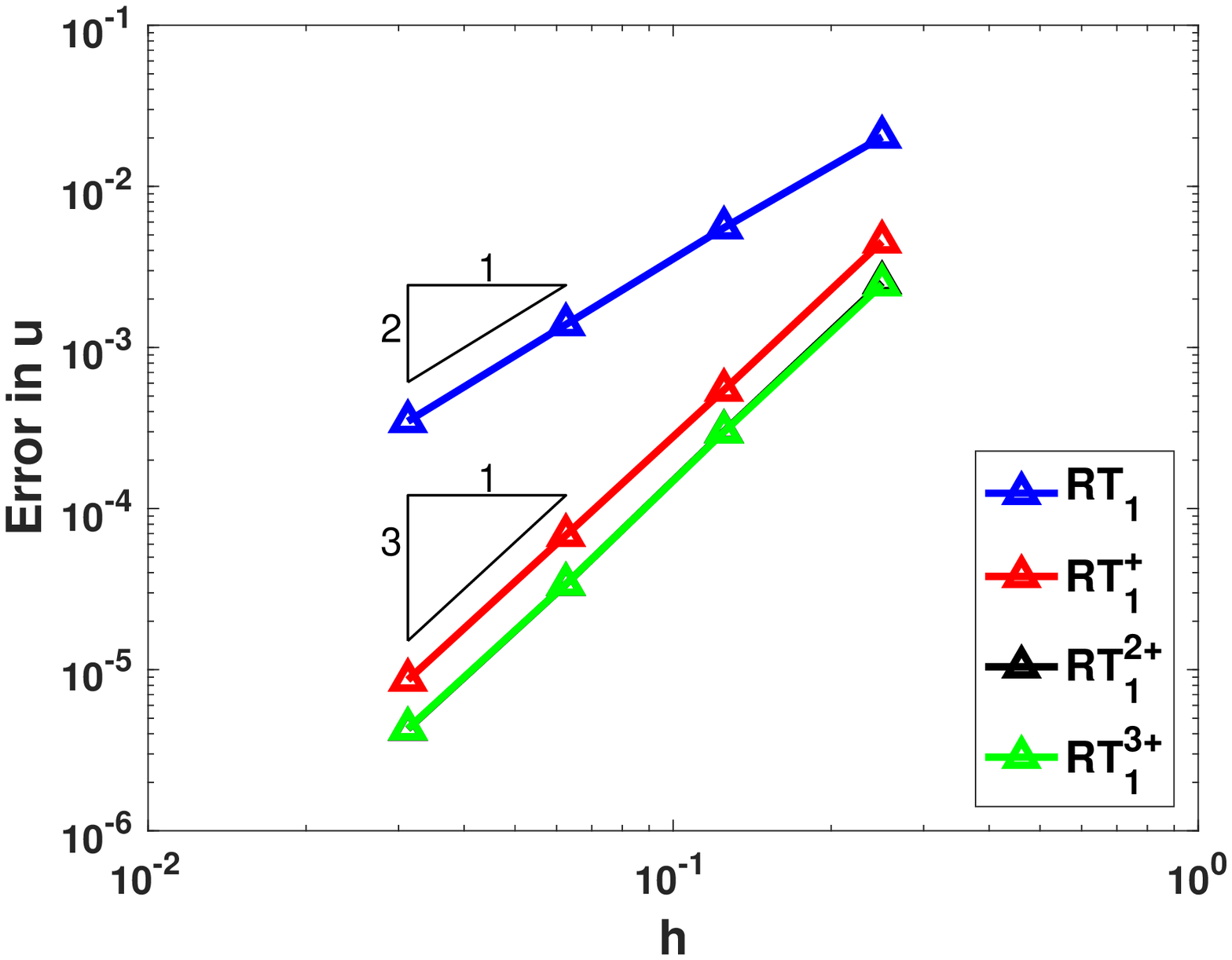} & \includegraphics[scale=0.35]{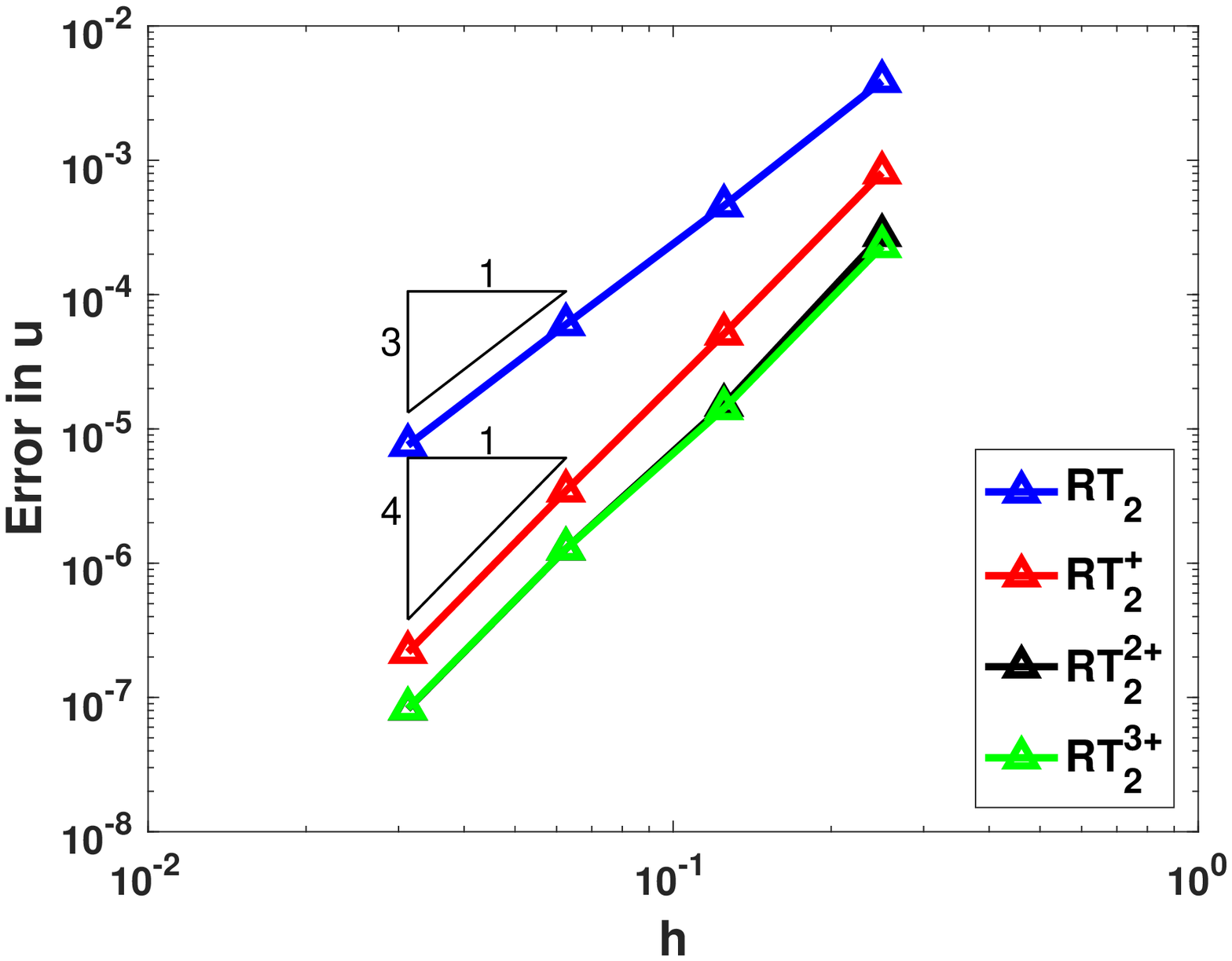}\tabularnewline
\includegraphics[scale=0.35]{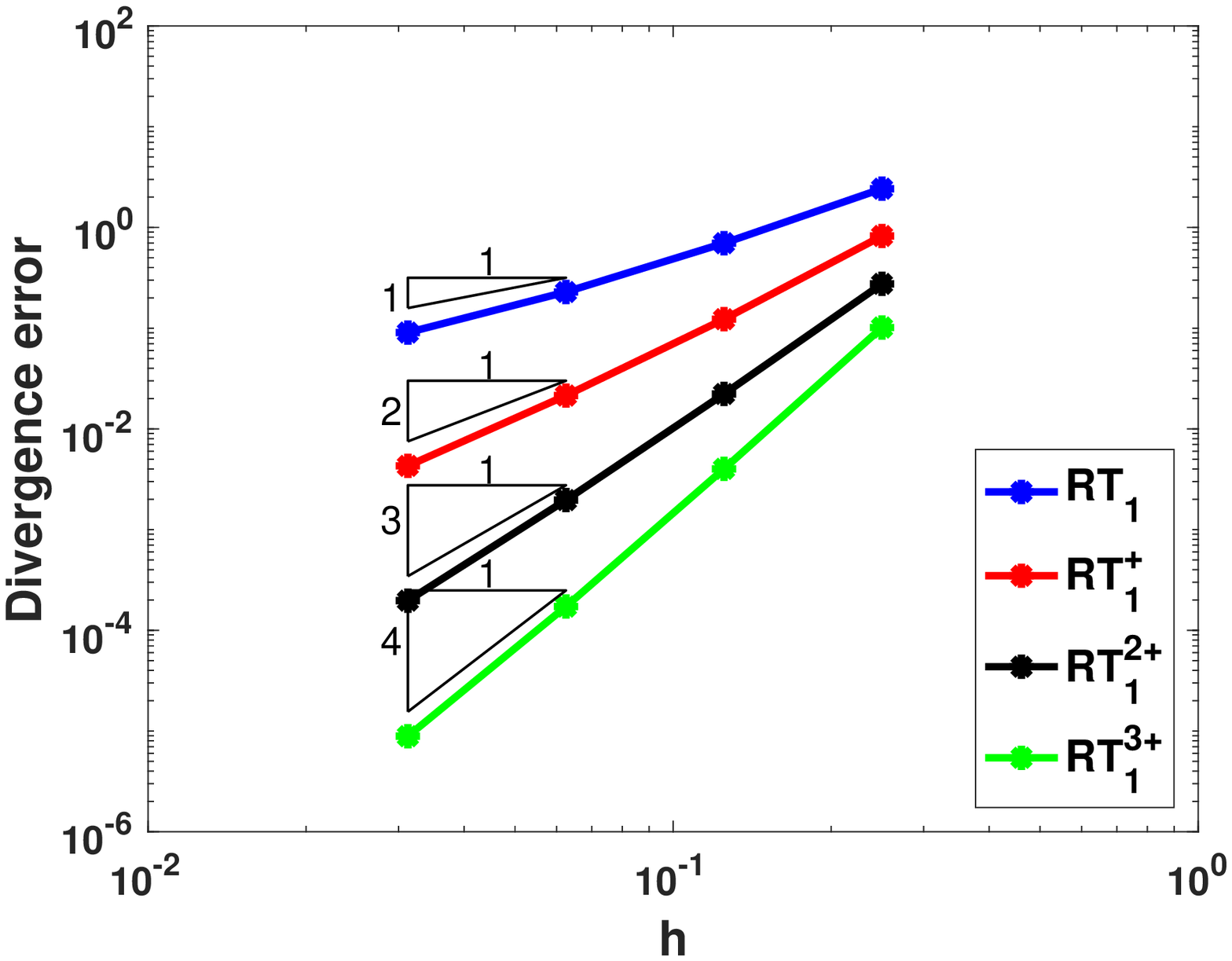} & \includegraphics[scale=0.35]{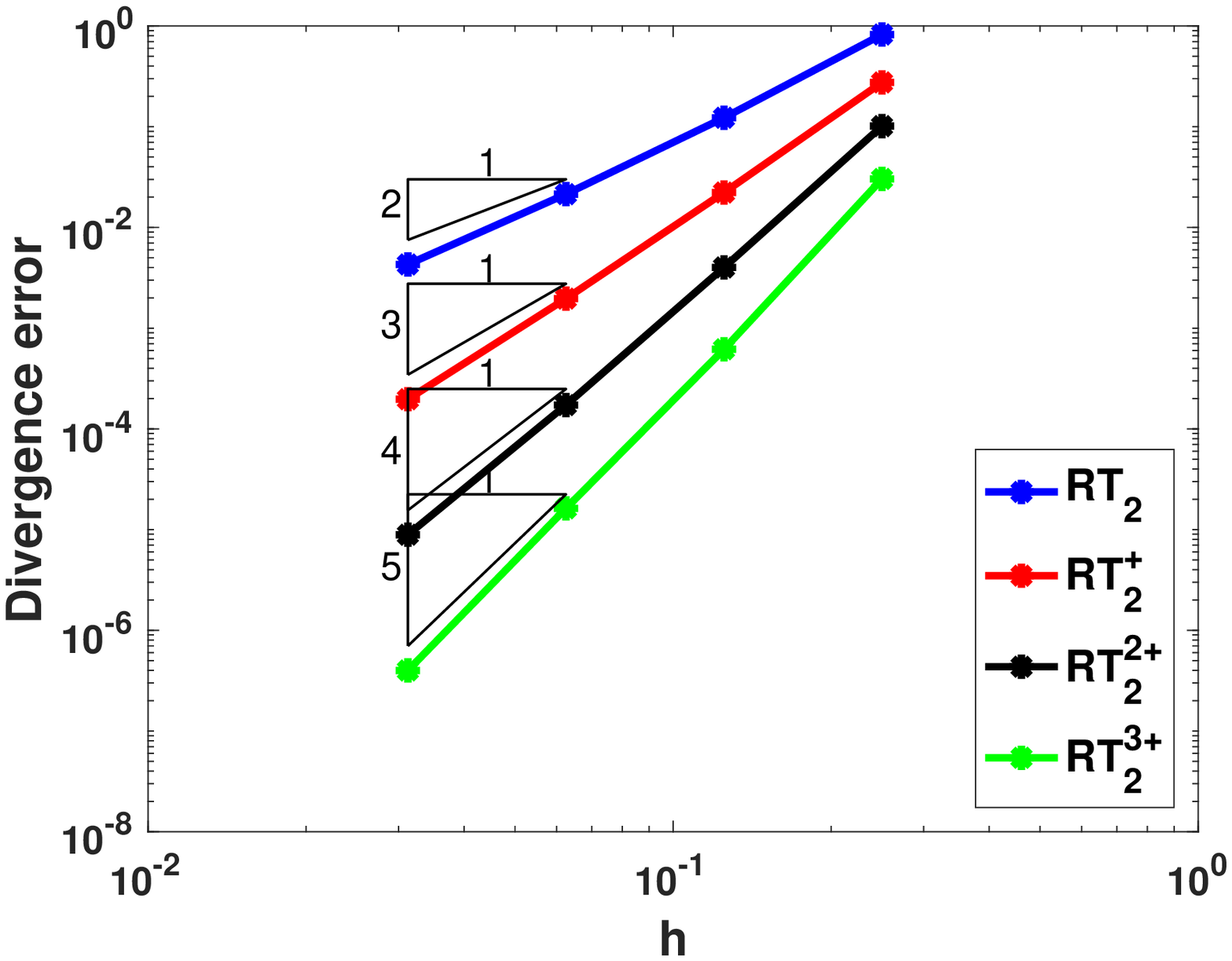}\tabularnewline
\end{tabular}
\par\end{centering}
\caption{$L^2$-errors in the flux (top), in $u$ (middle) and divergence (bottom)
with $k=$1 (left) and $k=2$ (right), using trapezoidal meshes and approximation space configurations
$RT_k^{n+}$, for $n=1,\cdots, 3$.\label{fig:Errors-Trapk12}}
\end{figure}

\begin{figure}
\begin{centering}
\begin{tabular}{cc}
\includegraphics[scale=0.35]{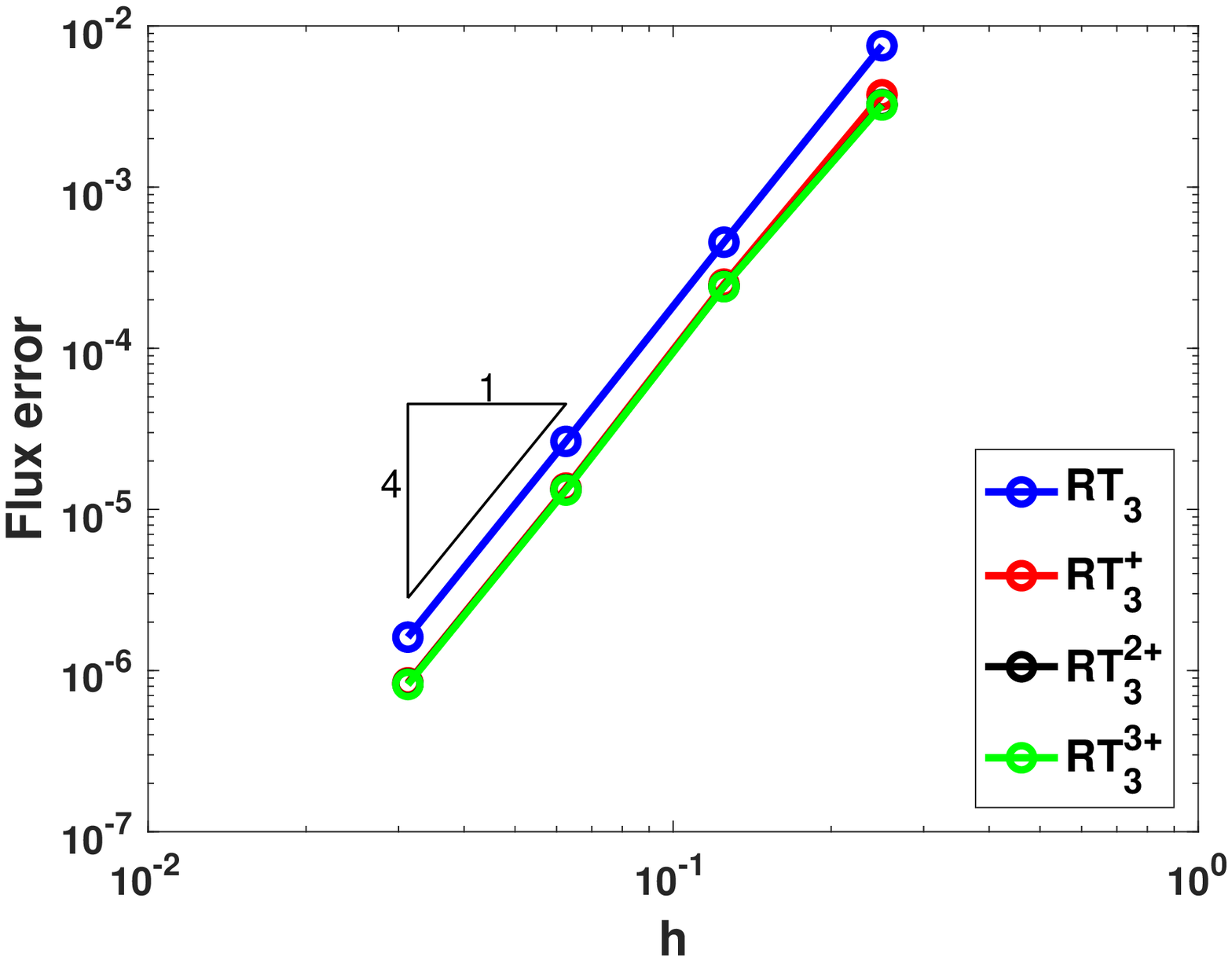} & \includegraphics[scale=0.35]{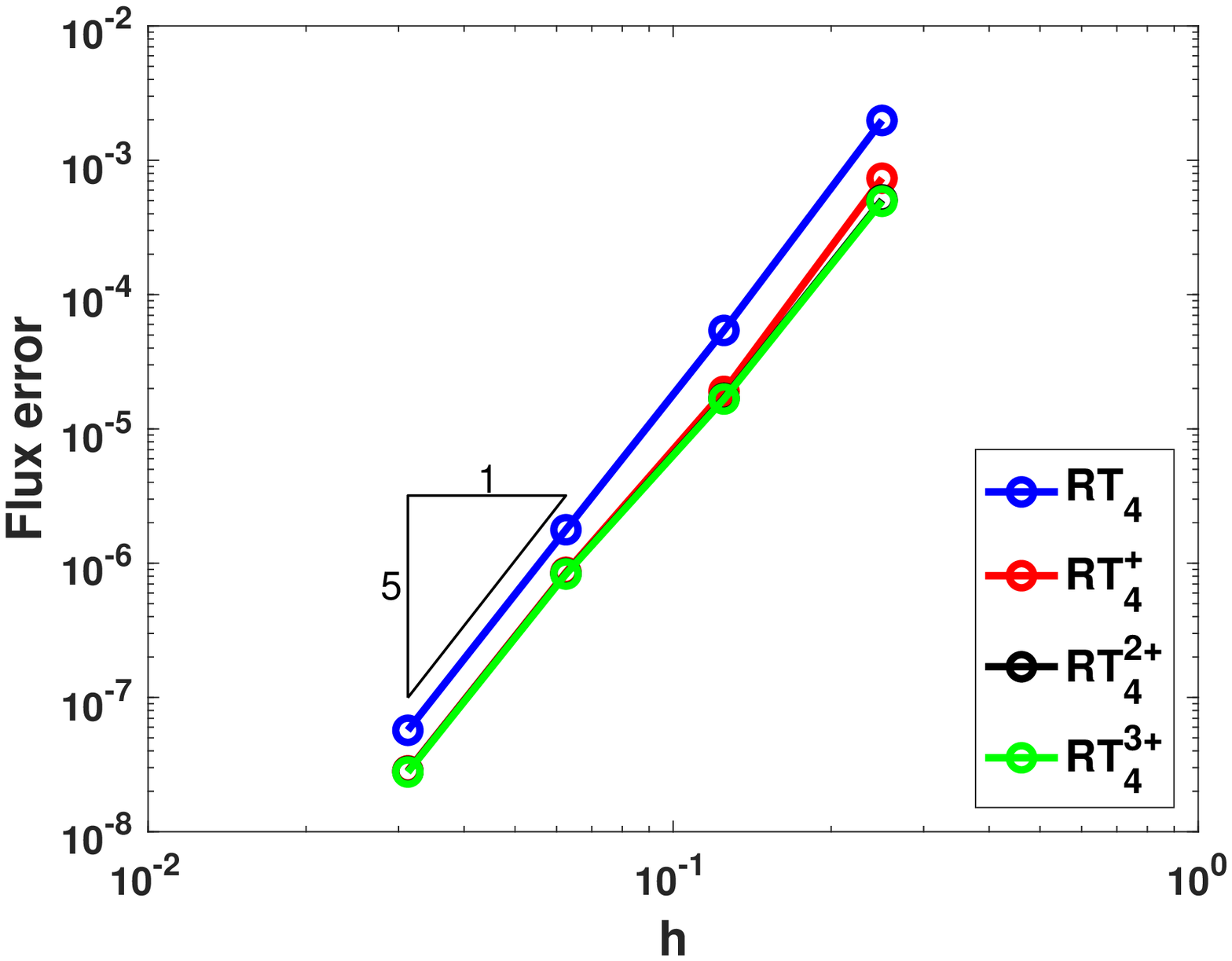}\tabularnewline
\includegraphics[scale=0.35]{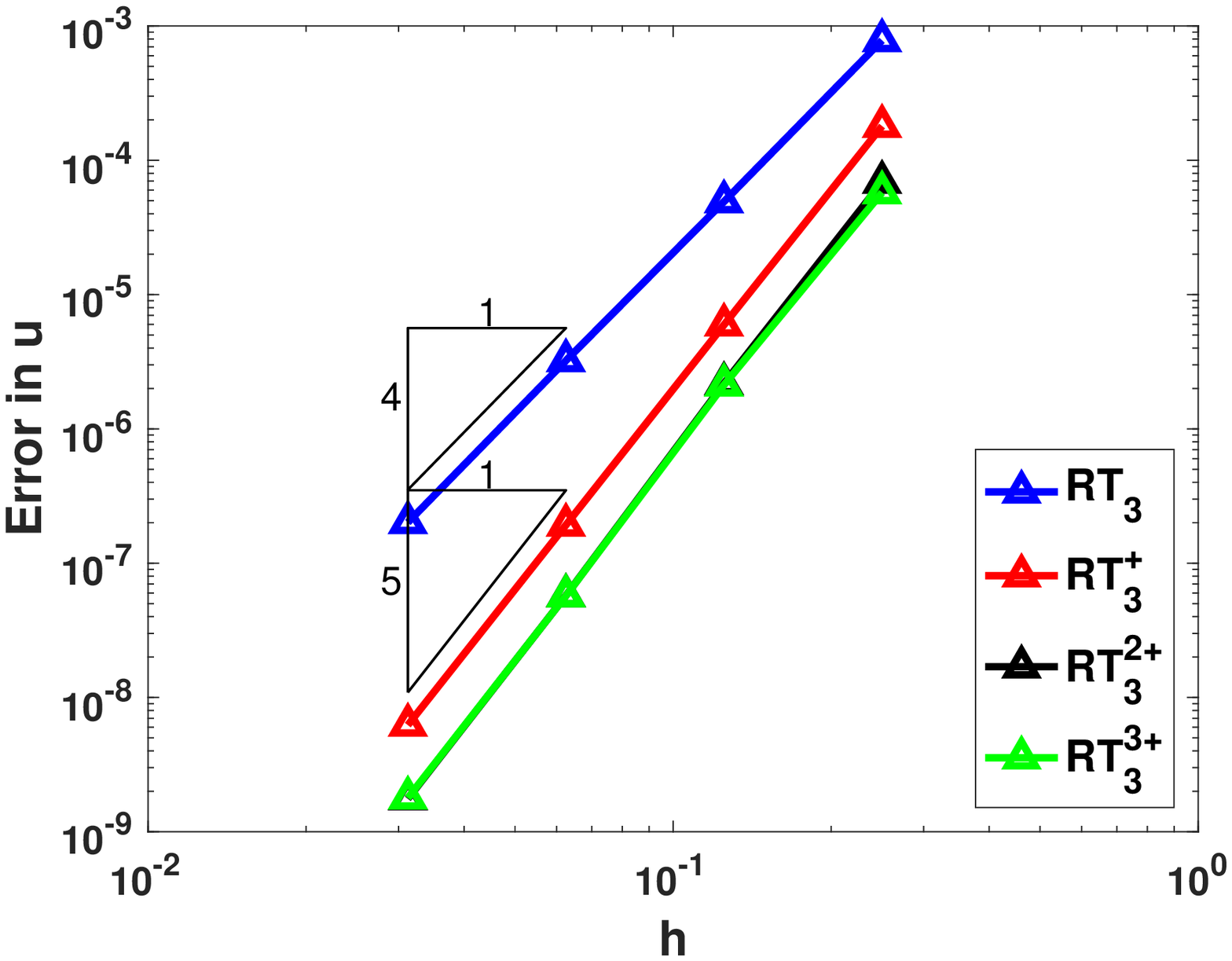} & \includegraphics[scale=0.35]{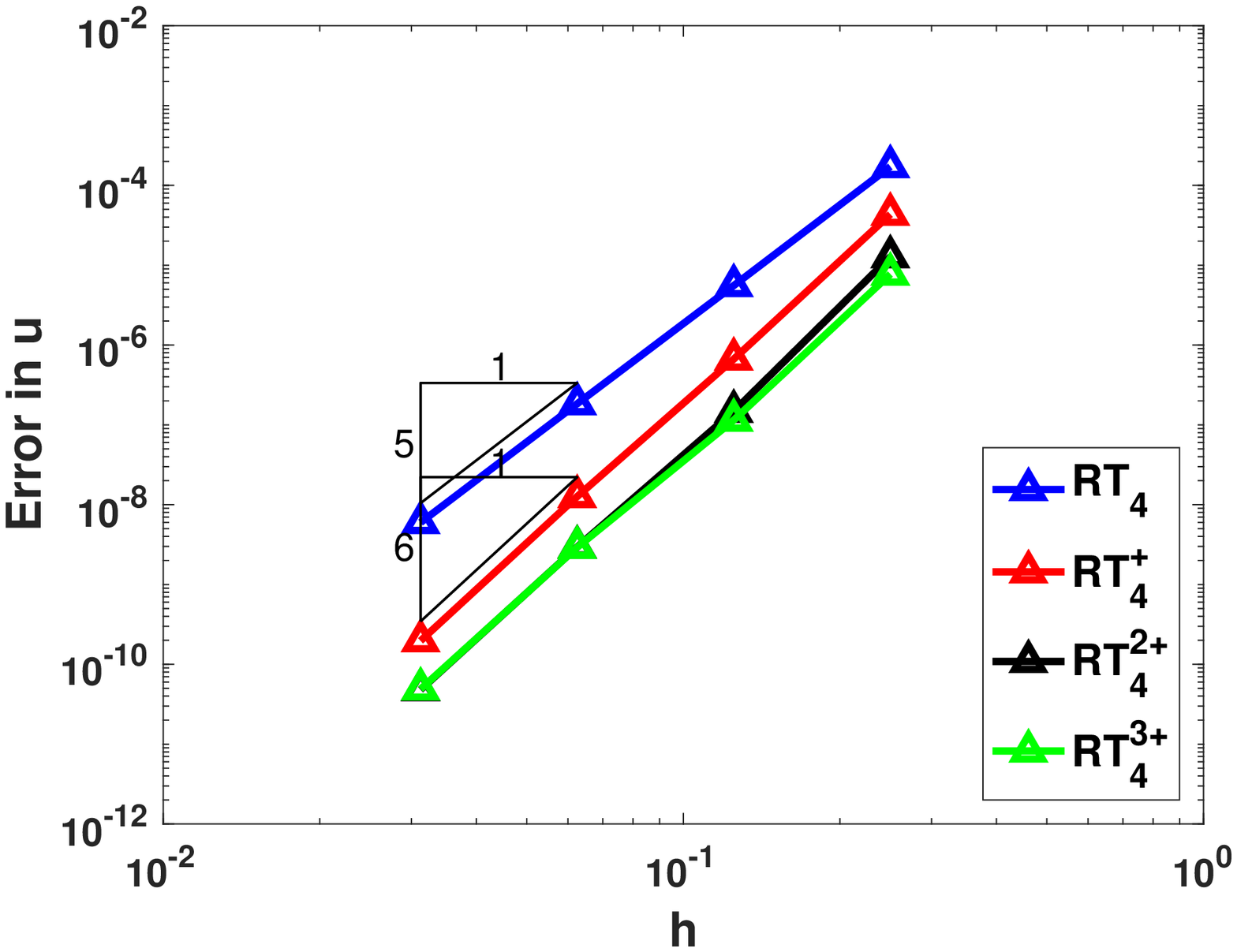}\tabularnewline
\includegraphics[scale=0.35]{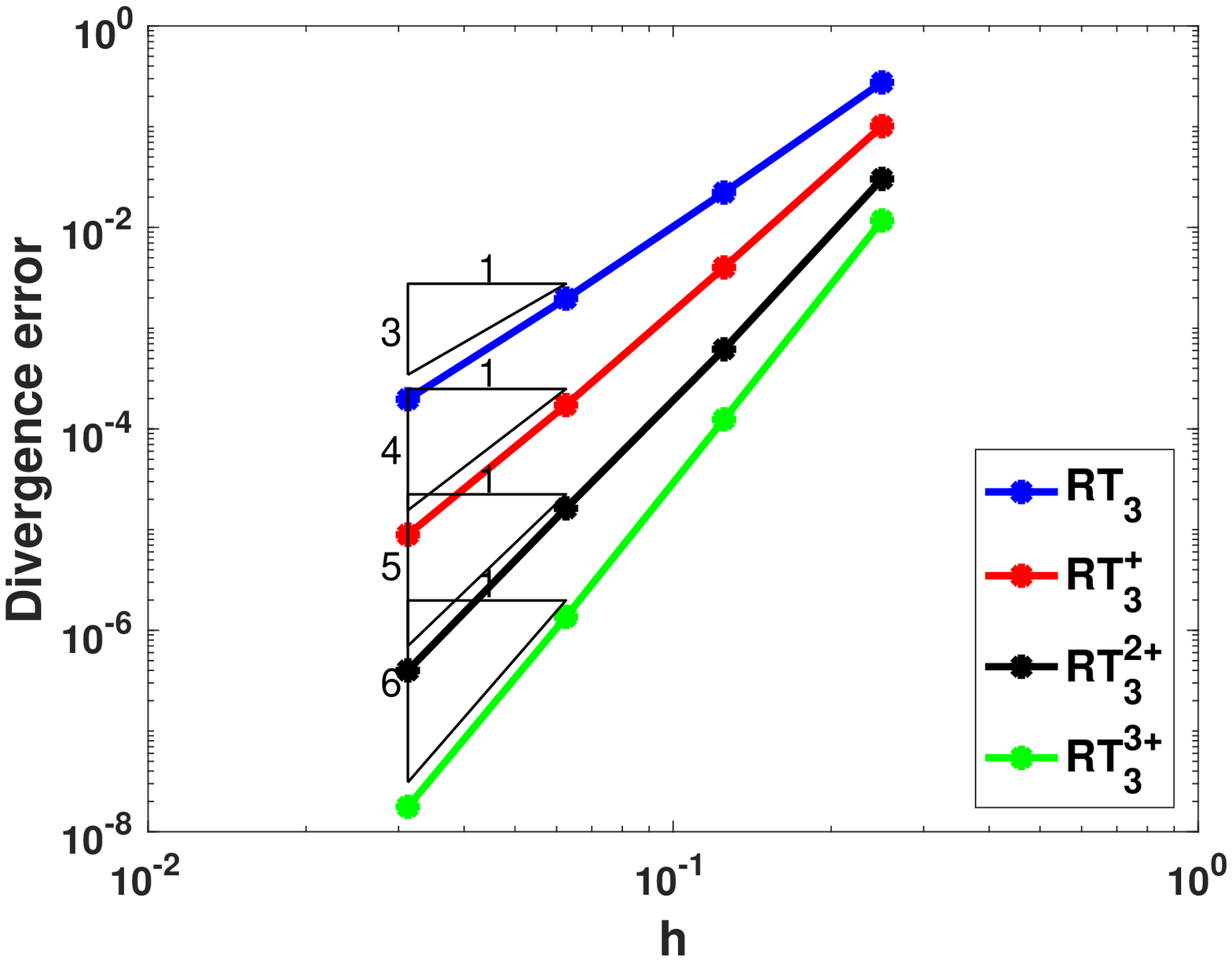} & \includegraphics[scale=0.35]{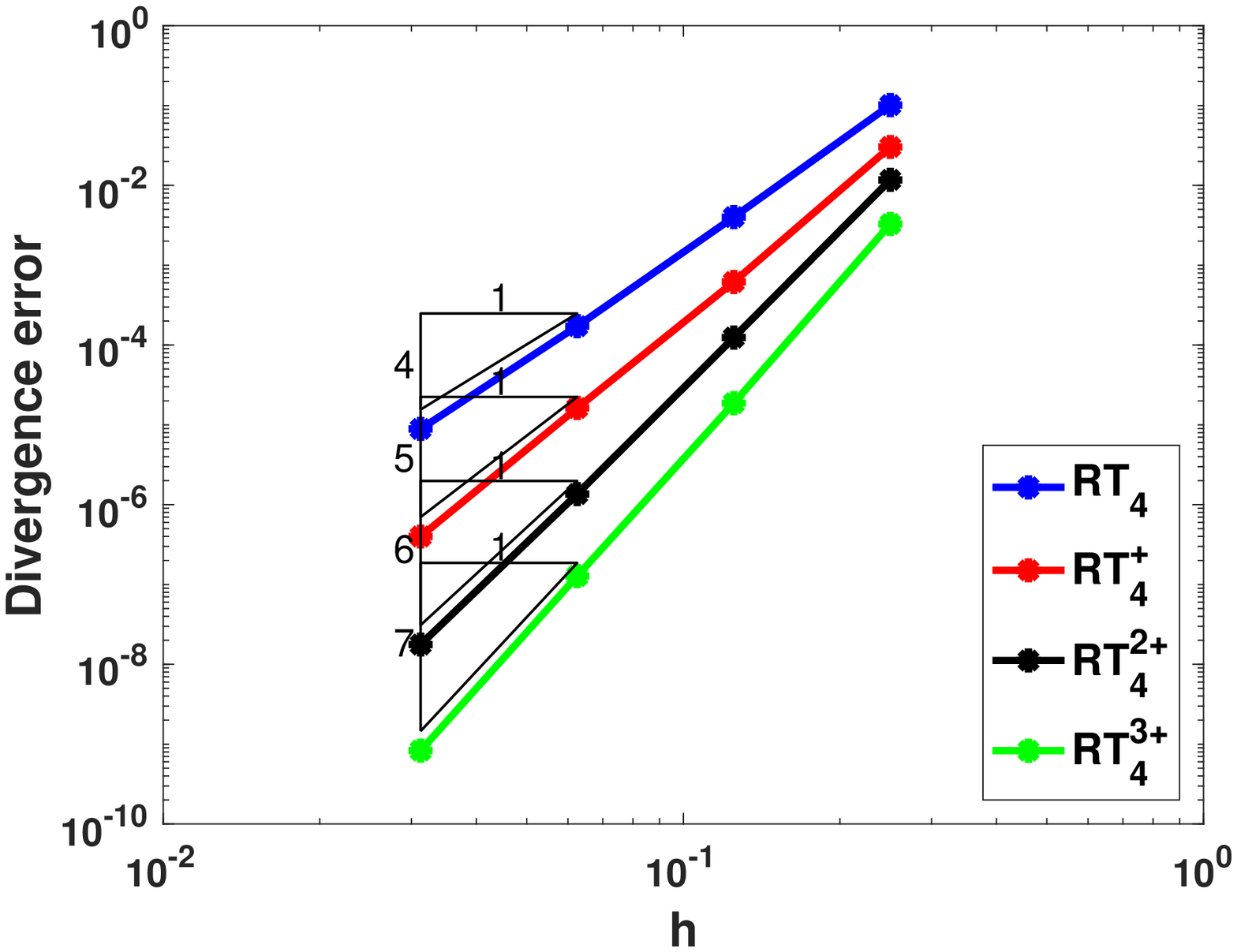}\tabularnewline
\end{tabular}
\par\end{centering}
\caption{$L^2$-errors in the flux (top), in $u$ (middle) and divergence (bottom)
with $k=$3 (left) and 4 (right), using trapezoidal meshes and approximation space configurations
$RT_k^{n+}$, for $n=1,\cdots ,3$.\label{fig:Errors-Trapk34}}
\end{figure}

The result  for  non-affine trapezoidal meshes are shown in Figure \ref{fig:Errors-Trapk12}, with $k=1$ and $k=2$, and in Figure \ref{fig:Errors-Trapk34}, for $k=3$ and $k=4$. Similar observations concerning flux and potential magnitude and accuracy orders hold as in the  affine cases. For the flux divergence, the effect of using non-affine meshes can be observed, by the reduction of one unit in the accuracy order to $k+n$ for the trapezoidal meshes, instead of $k+n+1$ in the affine case.

\section{Conclusions}\label{section6}

A general methodology is developed that is applicable to
 any space configuration $\hat{\mathbf{V}}_k, \hat{U}_k$,  defined in the master element $\hat{K}$, which is supposed to be used for the construction of divergence conforming  approximations for vector functions and for their divergence.  New insight on approximation properties of the divergence operator
 in mixed finite elements is obtained by showing that it is possible to  enrich the spaces   and get   discretizations of the  divergence operator as accurate as desired, without increasing the size of the condensed linear systems to be solved. 

 The main requirement is that the original 
 scheme  should allow  a direct factorization $\hat{\mathbf{V}}_k=\hat{\mathbf{V}}^\partial_k  \oplus \mathring{\hat{\mathbf{V}}}_k$  of  edge and internal flux functions, and that de Rham compatibility formula $\hat{U}_k = \nabla \cdot \hat{\mathbf{V}}_k$ holds. The enriched versions, 
 which results in higher order convergence rates of
 the divergence of the flux operator, are defined as
$ \hat{\mathbf{V}}_{k}^{n+} = \hat{\mathbf{V}}^{\partial}_{k}\oplus\mathring{\hat{\mathbf{V}}}_{k+n}$, with corresponding scalar spaces  $\hat{U}_k^{n+} =  \hat{U}_{k+n}$.  Projections allowing to demonstrate the commuting
property of the  de Rham diagram can be  extended to the enriched context.
From the coding point of view, the main requirement is that
the polynomial order of the edge flux functions 
can be chosen independently from the polynomial order
of the internal flux space. The enriched space  $ \hat{\mathbf{V}}_{k}^{n+}$ is simply obtained by pruning from $\hat{\mathbf{V}}_{k+n}$ the edge functions of degree $>k$.  For applications to mixed formulations, the computational cost of
matrix assembly increases with $n$, but the global condensed systems to be solved have
same dimension and structure of the original scheme.

\section*{Acknowledgments}
The authors  P. R. B. Devloo and S. M. Gomes thankfully acknowledges financial support  from FAPESP - Research Foundation of the State of S\~{a}o Paulo, Brazil (grant 2016/05144-0), and from CNPq - Brazilian Research Council (grants 305425/2013-7, 305823-2017-5, and 304029/2013-0, 306167/2017-4).


\section*{References}
\begin{thebibliography}{8}

\bibitem{ABF2002} D.N. Arnold, D. Boffi, and R.S. Falk, Approximation by quadrilateral finite elements,
Math. Comp., 71 (2002), pp. 909–922.

\bibitem{ABF}D.N. Arnold, D. Boffi, and R.S. Falk, Quadrilateral {\bf H}(\text{div}) finite elements,  SIAM J. Numer. Anal. 42 (6): 2429-2451, 2005.
%


\bibitem{BDFM}
F. Brezzi, J. Douglas, M. Fortin, L.D. Marini, Efficient rectangular mixed finite elements in two and three space variable, RAIRO Mod\'{e}l,  Math. Anal. Num\'{e}r. 2: 581-604, 1987.

\bibitem{BDM}F. Brezzi, J. Douglas and L.D. Marini, Two Families of Mixed Finite Elements for Second Order Elliptic Problems,  Numer. Math. 47: 217-235, 1985.


\bibitem{BrezziFortin1991} F. Brezzi, M. Fortin. 1991. \textit{ Mixed and Hybrid Finite Element Methods}, 
Springer Series in Computational Mathematics, vol. 15, Springer-Verlag, New York.


\bibitem{Surface} D.A. Castro, P.R.B. Devloo, A.M. Farias, S.M. Gomes, O. Dur\'an, Hierarchical high order finite element bases for $\mathbf{H}(\text{div})$ spaces based on
curved meshes for two-dimensional regions or manifolds, J. Comput. Appl. Math. 301: 241–258, 2016.

\bibitem{Hdiv3D} D. A.  Castro,  P. R. B. Devloo, A. M. Farias,  S. M. Gomes,  D. Siqueira, O.   Dur\'an. Three dimensional
hierarchical mixed finite element approximations with enhanced primal variable accuracy. Comput. Methods  Appl. Mech. and  Engrg.  306 (2016), 479--502. 

\bibitem{Jay2005}B. Cockburn, J. Gopalakrishnan, Error analysis of
variable degree mixed methods for elliptic problems via hybridization.
Mathematics of Computation 74 (252) (2005),  1653-1677.

\bibitem{DemkowiczHP}L. Demkowicz, P. Monk, L.Vardapetvan, W. Rachowicz.
De Rham diagram for $hp$ finite element spaces. Comput. Math. Appl.,
39 (2000), 29-38.

\bibitem{Hdiv2D2017}A.M. Farias, P.R.B. Devloo,  S.M.
Gomes, D. de Siqueira, D. A. Castro, Two dimensional mixed finite element approximations for elliptic problems with enhanced accuracy for  the potential and  flux divergence. Comput. Math.  Appl.  74: 3283–3295, 2017.

\bibitem{NeoPZ}http://github.com/labmec/neopz

\bibitem{RT} 
P.A. Raviart and J.M. Thomas, A mixed finite element method for 2nd order elliptic problems.  Lect. Not. Math. 606: 292-315, 1997.

\bibitem{Siqueira}
D. Siqueira, P.R.B. Devloo, S.M. Gomes. A new procedure for the construction of hierarchical high order Hdiv and Hcurl finite element spaces,  Journal of Computational and Applied Mathematics, 240: 204-214, 2013.

\end{thebibliography}
\end{document}